\documentclass[a4paper]{article}
\usepackage{amsthm,amsmath,amssymb}
\newtheorem{teor}{Theorem}[section]
\newtheorem{defin}[teor]{Definition}
\newtheorem{lemm}[teor]{Lemma}
\newtheorem{osse}[teor]{Remark}
\newtheorem{prop}[teor]{Proposition}
\newtheorem{defi}[teor]{Definition}
\newtheorem{coro}[teor]{Corollary}
\newtheorem{prob}[teor]{Problem}

\newcommand{\bele}{\begin{lemm}\begin{sl}}
\newcommand{\enle}{\end{sl}\end{lemm}}
\newcommand{\bedef}{\begin{defi}\begin{sl}}
\newcommand{\eddef}{\end{sl}\end{defi}}
\newcommand{\bete}{\begin{teor}\begin{sl}}
\newcommand{\ente}{\end{sl}\end{teor}}
\newcommand{\beos}{\begin{osse}\begin{rm}}
\newcommand{\eddos}{\end{rm}\end{osse}}
\newcommand{\bepr}{\begin{prop}\begin{sl}}
\newcommand{\empr}{\end{sl}\end{prop}}
\newcommand{\bepro}{\begin{prob}\begin{rm}}
\newcommand{\empro}{\end{rm}\end{prob}}
\newcommand{\bede}{\begin{defin}\begin{sl}}
\newcommand{\edde}{\end{sl}\end{defin}}
\newcommand{\beco}{\begin{coro}\begin{sl}}
\newcommand{\enco}{\end{sl}\end{coro}}

\newcommand{\quand}{\quad\text{and}\quad}

\newcommand{\quext}{\quad\text}

\newcommand{\RR}{\mathbb{R}}
\newcommand{\NN}{\mathbb{N}}

\newcommand{\beeq}[1]{\begin{equation}\label{#1}}
\newcommand{\eddeq}{\end{equation}}

\newcommand{\beeqa}[1]{\begin{eqnarray}\label{#1}}
\newcommand{\eddeqa}{\end{eqnarray}}

\newcommand{\beal}[1]{\begin{align}\label{#1}}
\newcommand{\eddal}{\end{align}}

\newcommand{\bespl}[1]{\begin{split}\label{#1}}
\newcommand{\edspl}{\end{split}}

\newcommand{\bega}[1]{\begin{gather}\label{#1}}
\newcommand{\edga}{\end{gather}}

\newcommand{\beeqax}{\begin{eqnarray*}}
\newcommand{\eddeqax}{\end{eqnarray*}}

\def\qed{\ifmmode 
  \else \leavevmode\unskip\penalty9999 \hbox{}\nobreak\hfill
  \fi
  \quad\hbox{\hskip.5em\vrule width.4em height.6em depth.05em\hskip.1em}}
\def\endproofsym{\qed}
\newcommand{\dimbox}{\hbox{\hskip.5em\vrule width.4em height.6em depth.05em\hskip.1em}}
\renewenvironment{proof}[1][Proof]{\trivlist\item[\hskip\labelsep{\hskip0pt
    {\normalfont\scshape#1.}\hskip .321429\parindent}]\ignorespaces}
{\endproofsym\endtrivlist}
\def\endnobox{\def\endproofsym{}\end{proof}\def\endproofsym{\qed}}

\newcommand{\no}{\nonumber}

\newcommand{\beeqao}{\begin{eqnarray}\no}
\newcommand{\bealo}{\begin{align}\no}
\newcommand{\besplo}{\begin{split}\no}
\newcommand{\begao}{\begin{gather}\no}

\newcommand{\duav}[1]{\langle{#1}\rangle}

\newcommand{\duavg}[1]{\left\langle{#1}\right\rangle}

\newcommand{\cc}{{\mathfrak c}}

\newcommand{\dtt}{_{tt}}

\newcommand{\perogni}{\forall\,}
\newcommand{\esiste}{\exists\,}
\newcommand{\itt}{\int_0^t}

\newcommand{\io}{\int_\Omega}

\newcommand{\OO}{_{\Omega}}

\newcommand{\lhs}{left hand side}
\newcommand{\rhs}{right hand side}

\DeclareMathOperator{\dive}{div}
\DeclareMathOperator{\deriv}{d}

\DeclareMathOperator{\dist}{dist}

\DeclareMathOperator{\spa}{span}

\DeclareMathOperator{\loc}{loc}

\DeclareMathOperator{\discr}{discr}

\newcommand{\CZH}{C^0([0,T];H)}

\let\TeXchi\chi
\def\chi{{\setbox0 \hbox{\mathsurround0pt
$\TeXchi$}\hbox{\raise\dp0 \copy0 }}}

\newcommand{\znn}{_{0,n}}
\newcommand{\unn}{_{1,n}}
\newcommand{\ii}{_\infty}

\newcommand{\nnt}{_{n,t}}

\newcommand{\zzn}{_{0,n}}

\newcommand{\calH}{{\mathcal H}}

\newcommand{\calG}{{\mathcal G}}

\newcommand{\calA}{{\mathcal A}}
\newcommand{\calF}{{\mathcal F}}
\newcommand{\calE}{{\mathcal E}}

\newcommand{\calM}{{\mathcal M}}

\newcommand{\calB}{{\mathcal B}}

\newcommand{\calC}{{\mathcal C}}
\newcommand{\calS}{{\mathcal S}}

\newcommand{\calQ}{{\mathcal Q}}
\newcommand{\calV}{{\mathcal V}}
\newcommand{\calY}{{\mathcal Y}}
\newcommand{\calI}{{\mathcal I}}

\newcommand{\calZ}{{\mathcal Z}}
\newcommand{\calW}{{\mathcal W}}

\newcommand{\zzu}{_{0,1}}
\newcommand{\zzd}{_{0,2}}

\newcommand{\dit}{\deriv\!t}
\newcommand{\dis}{\deriv\!s}
\newcommand{\ditau}{\deriv\!\tau}

\newcommand{\dir}{\deriv\!r}

\newcommand{\ddt}{\frac{\deriv\!{}}{\dit}}

\newcommand{\uuu}{u_{1,1}}
\newcommand{\uud}{u_{1,2}}

\newcommand{\calL}{{\cal L}}

\numberwithin{equation}{section}
\begin{document}

\title{On the 2D Cahn-Hilliard equation with inertial term}

\author{Maurizio Grasselli\\
Dipartimento di Matematica, Politecnico di Milano,\\
Via E.~Bonardi~9, I-20133 Milano, Italy\\
E-mail: {\tt maurizio.grasselli@polimi.it}
\and
Giulio Schimperna\\
Dipartimento di Matematica, Universit\`a di Pavia,\\
Via Ferrata~1, I-27100 Pavia, Italy\\
E-mail: {\tt giusch04@unipv.it}\\
\and
Sergey Zelik\\
Department of Mathematics, University of Surrey,\\
Guildford, GU2 7XH, United Kingdom\\
E-mail: {\tt S.Zelik@surrey.ac.uk}
}

\maketitle
\begin{abstract}
 P.~Galenko et al. proposed a modified Cahn-Hilliard equation
 to model rapid spinodal decomposition in non-equilibrium phase separation processes. This equation
 contains an inertial term which causes the loss of any regularizing effect
 on the solutions. Here we consider an initial and boundary value problem for this equation
 in a two-dimensional bounded domain. We prove a number of results related to well-posedness
 and large time behavior of solutions.
 In particular, we analyze the existence of bounded absorbing sets in two
 different phase spaces and, correspondingly, we establish the existence of the global attractor.
 We also demonstrate the existence of an exponential attractor.
\end{abstract}


{\bf AMS (MOS) subject clas\-si\-fi\-ca\-tion:}~~35B40,
 35B41, 35Q99, 82C26


\section{Introduction}
\label{secintro}
The celebrated Cahn-Hilliard equation was proposed to describe phase
separation phenomena in binary systems \cite{CH}. The standard version reads
\beeq{CH0}
  u_t -\Delta(-\Delta u +f(u))=0,
\end{equation}
where $u$ represents the relative concentration of one species in a given
domain $\Omega\subset\RR^N$, $N\leq 3$, while $f$
is the derivative of a non-convex potential accounting for the presence
of two species (e.g., $f(u)=u(u^2-1)$). Many papers are devoted
to the mathematical analysis of this equation and many
different features have already been carefully analyzed. Here we recall,
in particular, well-posedness for different boundary conditions and/or singular
potentials, asymptotic behavior of solutions, existence
of global and exponential attractors, analysis of the stationary states.
We confine ourselves to quote only some contributions, namely,
\cite{BF,Dl,El,EZ,KNP,LZ,MW,MMW,MZ,MZ2,NCS,NC,NC2,NC3,PRZ,RZ,RH,vW,WW1,WW2,WW3,WZ}.

Among the phase transformations involved in phase separation, a
peculiar one is named spinodal decomposition, which indicates a
situation in which both the phases have an equivalent symmetry
and differ only in composition (see \cite{Ca}, cf.~also
\cite{Gr,MW}). It has been noted that, in certain materials
like glasses, \eqref{CH0} needs to be modified in order to
describe strongly non-equilibrium decomposition generated by
deep supercooling into the spinodal region (cf.
\cite{GL1,GL2,GL3} and references therein). In this respect,
P.~Galenko et al.\ proposed a modification based on the
relaxation of the diffusion flux (see, for instance,
\cite{Gal,GJ}) which yields the following evolution equation
\beeq{CHmod}
  \epsilon u\dtt + u_t -\Delta(-\Delta u +f(u))=0,
\end{equation}
where $\epsilon>0$ is a relaxation time. This equation shows a good agreement with experimental data
extracted from light scattering on spinodally decomposed glasses (cf. \cite{GL1}).

>From the mathematical viewpoint, equation \eqref{CHmod} was
firstly analyzed in \cite{D} as a dissipative dynamical system.
Those pioneering results were then improved and generalized in
\cite{ZM1,ZM2} with some restrictions on $\epsilon$. Then, in
\cite{GGMP1}, the existence of a family of exponential
attractors, robust as $\epsilon$ goes to $0$, was established.
All the quoted papers were devoted to the one-dimensional case
which is relatively easy since a rather weak (energy
bounded) solution is also bounded in $L^\infty$. This is no longer true in
dimension two or three and, while the existence of an energy
bounded solution can be still proven rather easily (with some
restrictions on $f$), its uniqueness or the existence of
smoother solutions appear nontrivial.

Equation \eqref{CHmod} behaves in a very different
way with respect to \eqref{CH0} since there is
no regularization of the solutions in finite time. This smoothing property
can be restored if one adds a viscosity term of the form $-\alpha \Delta u_t$,
with $\alpha>0$ (see \cite{NC} for a physical justification). The viscous variant
of \eqref{CHmod} in dimension three was firstly analyzed in detail in \cite{GGMP2}.
In particular,
the authors constructed a family of exponential attractors which is robust as $(\alpha,\epsilon)$
goes to $(0,0),$ provided that $\epsilon$ is dominated by $\alpha$. For other recent
contributions related to the viscous version of \eqref{CHmod}, the reader is referred
to \cite{Bo,BGM,GPS,Ka}.
Going back to \eqref{CHmod}, in the three-dimensional case,
an analysis of the longtime behavior of the solutions based only
on the existence result was carried out in \cite{S}.
Then, a nonisothermal phase-separation system with memory was considered in \cite{V}
(see also \cite{LR} for existence of weak solutions).
This system can be easily reduced to \eqref{CHmod}
by neglecting the temperature effects
and taking the memory kernel
equal to a decreasing exponential in the remaining equation.
However, the assumptions made on the memory kernel exclude this possibility since
the decreasing exponential is not $\theta$-sectorial with $\theta\in(0,\frac{\pi}{2})$
(see \cite[Def.~1 and assumption (P0)]{V}). We recall that this hypothesis is crucial to ensure
the parabolic nature of the integrodifferential equation.

Hence, many questions about \eqref{CHmod} are still unanswered so far, namely, uniqueness
of weak solutions, well-posedness in stronger settings, construction
of bounded absorbing sets, existence of global attractors as well as exponential attractors.
The present contribution gives several answers in the
two-dimensional case. Here we are not concerned with the dependence on $\epsilon$,
thus we take $\epsilon=1$. This dependence will be possibly studied in a future paper.
Regarding the three-dimensional case we observe that the smallness of $\epsilon$
seems to play a crucial role if we want to extend the mentioned results (see \cite{GSSZ})
as in the case of damped wave equation with supercritical nonlinearities treated in \cite{Ze}.
In two dimensions we can avoid the restriction on $\epsilon$ since we can take advantage of
the Br\'ezis-Gallouet inequality (cf. \cite{BG}).
%
%
It is worth noting that this restricts our analysis to
functions $f$ with (at most) a cubic growth at infinity. In this
sense, the situation we meet for the 2D fourth order equation
\eqref{CH0} is surprisingly similar to what happens for
the 3D (second order) damped wave equation studied, e.g.,
in \cite{ACH,Ba2,GP,PZ}, where
an analogous growth restriction on $f$ is assumed.
%
%
Therefore, in this paper we will study the
initial and boundary value problem
\begin{align}
  \label{CHin}
  &u\dtt(t)+u_t(t)-\Delta(-\Delta u(t)+f(u(t)))=g,\qquad\textrm{ in }\, \Omega, \;t> 0,\\
  \label{bc}
  &u(t)=\Delta u(t)=0, \qquad\textrm{ on }\, \partial\Omega,\; t>0,\\
  \label{ic}
  &u(0)=u_0,\;u_t(0)=u_1, \qquad\textrm{ in }\, \Omega,
\end{align}
where $\Omega\subset\RR^2$ is a given smooth and bounded
domain. Here $f$ is the derivative of
a nonconvex smooth potential, while $g$ is a known time-independent source term.
The choice of the boundary conditions is somewhat artificial but it allows us
to simplify the presentation (see also \cite{GGMP1,ZM1,ZM2}). However, our arguments could be
recasted also when we deal with usual no-flux boundary conditions (like, e.g., in \cite{D,GPS}).

The paper is organized as follows. The next Section~\ref{secmain}
is devoted to the statement of our hypotheses and the proof
of well-posedness with initial data in $H^3\times H^1$ (quasi-strong solutions).
In Section~\ref{secV2}, dissipativity and existence of the global
attractor are demonstrated for the semiflow generated by the
solutions we previously found.
Then, a regularity property of the attractor is established in Section~\ref{secV2b}
while the existence of an exponential attractor is proven in Section~\ref{secexpo}.
The class of energy bounded solutions is carefully analyzed in Section~\ref{secV0},
showing, in particular, a well-posedness result and the existence of the global attractor.


\section{Existence and uniqueness of quasi-strong solutions}
\label{secmain}

Let us set $H:=L^{2}(\Omega)$ and denote by
$(\cdot,\cdot)$ the scalar product both in $H$ and
in $H\times H$, and by $\|\cdot\|$ the induced norm.
The symbol $\|\cdot\|_{X}$ will indicate the norm
in the generic real Banach space $X$.
Next, we set $V:=H^1_0(\Omega)$, so that
$V'=H^{-1}(\Omega)$ is the topological
dual of $V$. The space $V$ is endowed with
the scalar product
\beeq{defiBu}
  (\!(v,z)\!):=\io\nabla v\cdot\nabla z,
\end{equation}
and the related norm. We also denote by $A:D(A)\to H$
the Laplace operator with homogeneous Dirichlet boundary
condition. It is well known that $A$ is a strictly positive
operator with $D(A)=H^2(\Omega)\cap V$ (note that
we shall always suppose $\Omega$ smooth enough),
so that we can define, for $s\in\RR$, its powers
$A^s:D(A^s)\to H$. Moreover, we introduce the
scale of Hilbert spaces
\beeq{deficalVs}
  \calV_s:=D(A^{\frac{s+1}2})\times D(A^{\frac{s-1}2}),
\end{equation}
so that we have, in particular, $\calV_0=V\times V'$
and $\calV_1=(H^2(\Omega)\cap V)\times H$. The spaces $\calV_s$ are
naturally endowed with the graph norm
\beeq{definorVs}
  \|(u,v)\|_s^2 := \|A^{\frac{s+1}2} u\|^2
   +\|A^{\frac{s-1}2}v\|^2.
\end{equation}
Our hypotheses on the nonlinear function $f$ are the following:
%
\beal{f1}
  & f\in C^{2,1}_{\loc}(\RR;\RR), \quad f(0)=0, \quad
   \esiste r_0\ge0:~\, f(r)r \ge 0~~\perogni |r|\ge r_0,\\
 \label{f2}
  & \esiste \lambda\ge0
   :~~f'(r)\ge
   -\lambda~~\perogni r\in\RR,\\
 \label{f3}
  & \esiste M\ge 0:~~|f''(r)|\le M(1+|r|)~~\perogni r\in\RR.
\end{align}
We note by $F$ the primitive of $f$ such that $F(0)=0$.
%
%
Then, for any fixed final time $T>0$,
problem \eqref{CHin}-\eqref{ic}, noted in the sequel
as Problem~(P), can be written in the more abstract form
\begin{align}\label{CH}
  & u\dtt+u_t+A(A u+f(u))=g, \quext{in }\/V',~~\text{a.e.~in }(0,T),\\
 \label{iniz}
  & u|_{t=0}=u_0, \quad u_t|_{t=0}=u_1, \quext{a.e.~in }\,\Omega,
\end{align}
where $u_0$, $u_1$ are given initial data.
Formally testing \eqref{CH} by $A^{-1} u_t$, one readily sees
that the {\sl energy}\/ functional
\beeq{defiE}
  \calE:\calV_0\to \RR, \qquad
   \calE(u,v):=\frac12\|(u,v)\|_0^2+\io F(u)
     -\duav{g,A^{-1}u}
\end{equation}
can be associated with \eqref{CH} and,
due to \eqref{f3}, $\calE$ is
finite for all $(u,v)\in \calV_0$, provided
that $g\in V'$ (cf.~\eqref{regog} below; however,
at this stage $g\in D(A^{-3/2})$ would suffice).
Moreover, by \eqref{f1}, $F$ is bounded from below (and $\calE$ as well).

\beos\label{weakf1}
The last condition in \eqref{f1} is assumed just to avoid
further technicalities. Indeed, it can be relaxed by taking
$$
\liminf_{\vert r\vert \to \infty} \frac{f(r)}{r} > -\lambda_1,
$$
where $\lambda_1$ is the first eigenvalue of $A$. Note that, in this case, 
we can choose $F$ such that
$$
F(r) \geq -\frac{\kappa}{2} r^2,
$$
for some $\kappa < \lambda_1$. Thus $\cal E$ is still bounded from below.
\eddos

In order to distinguish the solutions according to their
smoothness, we introduce the following terminology. Given some
$T>0$, a solution $(u,u_t) \in L^\infty(0,T;\calV_0)$ to
\eqref{CH}-\eqref{iniz} will be named {\sl energy solution}.
Instead, $(u,u_t) \in L^\infty(0,T;\calV_1)$ will be called
{\sl weak solution}. Note that energy solutions are weaker than
weak solutions. Speaking of smoother solutions, $(u,u_t) \in
L^\infty(0,T;\calV_2)$ will be called {\sl quasi-strong
solution} and $(u,u_t) \in L^\infty(0,T;\calV_3)$ will be a
{\sl strong solution}. In the latter case, $u$ satisfies
equation \eqref{CH} almost everywhere in $\Omega\times (0,T)$.

It seems natural to look first for energy solutions.
Nonetheless, we prefer to investigate the class of quasi-strong solutions
and then construct the energy solutions
by an approximation-limit argument. Actually, we will see
in Section~\ref{secV0} that some properties of energy solutions
(like, e.g., uniqueness and asymptotic behavior) are rather delicate to handle.
Regarding weak solutions, there are still some open questions (see Remark \ref{ancheV1} below).

Our first result states then the well-posedness of Problem~(P) in the class of quasi-strong solutions.
\bete\label{teoesiforte}
 Let us assume\/ \eqref{f1}-\eqref{f3} and
 \begin{align}\label{regog}
   & g\in V',\\
  \label{regou0forte}
   & (u_0,u_1)\in\calV_2.
 \end{align}
 Then, there exists\/ {\sl one and only one}\/ function
 \beeq{regouforte}
   u\in W^{2,\infty}(0,T;V')\cap W^{1,\infty}(0,T;V)
   \cap L^\infty(0,T;D(A^{3/2}))
 \end{equation}
 solving\/ {\rm Problem~(P)}.
\ente
Before proving the theorem, let us observe that,
being $T>0$ arbitrary, $u$ can be thought to be defined for
all times $t\in(0,\infty)$. Moreover, we remark
that a solution to~(P) will be indifferently
noted in the sequel either as $u$ or as a couple
$(u,u_t)$, the latter notation being preferred when
we want to emphasize the role of some $\calV_s$
as a phase space. We will also frequently write
$U$ for $(u,u_t)$. Moreover, throughout the remainder of
the paper the symbols $c$, $\kappa$, and $c_i$,
$i\in {\mathbb N}$, will denote positive constants
depending on the data $f,g$ of the problem, but independent
of the initial datum and of time. The value of $c$ and $\kappa$
is allowed to vary even within the same line.
Analogously, $Q:\RR^+\to\RR^+$ denotes a
generic monotone function.
Capital letters like $C$ or $C_i$ will be used to
indicate constant which have other dependencies
(in most cases, on the initial datum). Finally,
the symbol $c\OO$ will denote some embedding constants
depending only on the set $\Omega$.
%

\vspace{2mm}

\noindent%
{\bf Proof of Theorem~\ref{teoesiforte}.}~~%
To prove the regularity \eqref{regouforte},
we perform a number of a priori estimates. These
may have just a formal character in this setting, but could
be justified by working, e.g., in a Faedo-Galerkin
approximation scheme and then taking the limit.
The details of this standard procedure are omitted;
however, a sketch will be given in the proof
of Theorem~\ref{teoesidebo} below.

Thus, let us start with the energy estimate.
Let us set $U:=(u,u_t)$ and $U_0:=(u_0,u_1)$,
for brevity.
Testing \eqref{CH} by $A^{-1}u_t$, we get
\beeq{conto11}
  \ddt\calE(U)
   +\|u_t\|_{V'}^2\le 0.
\end{equation}
Then, let us take a (small) constant $\beta>0$, test
\eqref{CH} by $\beta A^{-1}u$, and add the result
to \eqref{conto11}. Noting that, by \eqref{f2},
$f(r)r\ge F(r)-\lambda r^2/2$, it is then
not difficult to infer the dissipativity of the
energy, i.e.,
\beeq{st11}
  \calE(U(t))
   \le \calE(U_0)e^{-\kappa t}
   + Q(\|g\|_{V'}).
\end{equation}
Hence, being by \eqref{f3}
$\io F(u_0)\le Q(\|u_0\|_V)$
and $\calE(u,v)\le Q(\|(u,v)\|_0)+Q(\|g\|_{V'})$,
and recalling \eqref{defiE}, we also derive
\beeq{st12}
  \|U(t)\|_0
   \le Q(\|U_0\|_0)e^{-\kappa t} + Q(\|g\|_{V'}).
\end{equation}
Next, by \eqref{conto11} we obtain that $\|u_t\|_{V'}$
is summable over $(0,\infty)$. More precisely,
integrating that relation from
an arbitrary $t\ge0$ to $\infty$,
and using \eqref{st11}-\eqref{st12}, we infer
\beeq{dissin}
  \int_t^\infty\|u_t(s)\|_{V'}^2\,\dis
   \le \calE(U(t))-\calE\ii
   \le Q(\|U_0\|_0)e^{-\kappa t} + Q(\|g\|_{V'}),
\end{equation}
where $\calE\ii$ is the limit for $t\nearrow\infty$
of the (nonincreasing) function
$t\mapsto\calE(U(t))$.
%
%
>From this moment on, let us denote, for brevity, by $\cc$ the
\rhs\ of \eqref{st12}.
%
%
Let us then differentiate Problem~(P) with respect to time. Again, this
formal procedure can be justified in the
Faedo-Galerkin approximation.
Setting $v:=u_t$, we obtain
\begin{align}\label{CH'}
  & v\dtt+v_t+A(Av+f'(u)v)=0,\\
 \label{iniz'1}
  & v|_{t=0}=v_0:=u_1,\\
 \label{iniz'2}
  & v_t|_{t=0}=v_1:=-u_1-A^2u_0-Af(u_0)+g.
\end{align}
Being $u_0\in D(A^{3/2})$ by \eqref{regou0forte}
and owing to \eqref{f1}, it is
not difficult to check that $v_1\in V'=D(A^{-1/2})$. More
precisely, we have that
\beeq{u0v0}
  \|V_0\|_0\le Q(\|U_0\|_2)
\end{equation}
(we have set here $V_0:=(v_0,v_1)$ and
$V:=(v,v_t)$).
Next, let us test \eqref{CH'} by $A^{-1}(v_t+\beta v)$,
with small $\beta>0$ as before. This gives
\bealo
  & \ddt\Big(\frac12\|V\|_0^2
      +\beta\duav{v_t,A^{-1}v}
      +\frac\beta2\|v\|_{V'}^2 \Big)
   +(1-\beta)\|v_t\|_{V'}^2\\
 \label{conto21}
  & \mbox{}~~~~~
   +\beta\|v\|_V^2
   +\big(f'(u)v,v_t\big)
   +\beta\big(f'(u)v,v\big)
   \le 0.
\end{align}
A straightforward computation then shows that
\beeq{conto22}
  \big(f'(u)v,v_t\big)
   =\frac12\ddt\big(f'(u)v,v\big)-
    \frac12\big(f''(u)v^2,v\big).
\end{equation}
Let us then now pick $\beta$ in \eqref{conto21}
so small that
\beeq{inconto21}
  \frac12\|V\|_0^2
   +\beta\duav{v_t,A^{-1}v}
   +\frac\beta2\|v\|_{V'}^2
  \ge\frac14\|V\|_0^2.
\end{equation}
Next, for $L>0$ whose value is to be chosen later,
we have by interpolation
\beeq{conto24.2}
  L\ddt\|v\|_{V'}^2+
   2\beta L\|v\|_{V'}^2 
  \le \frac\beta{8}\|V\|_0^2
   +c(L,\beta)\|v\|_{V'}^2
  \le \frac\beta{8}\|V\|_0^2
   +\cc.
\end{equation}
Then, let us add the above inequality to \eqref{conto21},
set
\beeq{deficalF}
  \calF:=\frac12\|V\|_0^2
    +\beta\duav{v_t,A^{-1}v}
    +\frac\beta2\|v\|_{V'}^2
    +\frac12\big(f'(u)v,v\big)
    +L\|v\|_{V'}^2,
\end{equation}
and observe that for a suitable choice of
$L$ (depending on $\lambda$ and
on $\beta$ taken before),
by interpolation there holds
\beeq{propcalF}
  \frac12\big(f'(u)v,v\big)
   +L\|v\|_{V'}^2
  \ge -\frac\lambda2\|v\|^2
   +L\|v\|_{V'}^2
  \ge -\frac18\|v\|_{V}^2.
\end{equation}
Hence, it is easy to see that $\calF$ satisfies, for
some $\sigma>0$ independent of the initial data,
%
%
\beeq{conto25}
  \calF\ge
   \sigma\|V\|_0^2.
\end{equation}
%
%
%
%
%
%
%
%
%
%
%
Moreover, recalling \eqref{conto22} and possibly
taking a smaller $\beta$, for some $\kappa>0$
independent of the initial data
we can rewrite \eqref{conto21} in the form
\beeq{conto21new}
  \ddt\calF
   +2\kappa\calF
  \le \cc+\frac12\big(f''(u)v^2,v\big).
\end{equation}
At this point, let us test \eqref{CH} by
$Au$. Using \eqref{f1}-\eqref{f3},
Sobolev embeddings and interpolation,
we can estimate the nonlinear term this way:
\bealo
  \duavg{Af(u),Au}
   & =\big(f'(u)\nabla u,\nabla\Delta u\big)
     \le \|f'(u)\nabla u\| \|u\|_{D(A^{3/2})}\\
 \no
   & \le c \big(1+\|u\|_{L^4(\Omega)}^2\big)
   \|\nabla u\|_{L^\infty(\Omega)} \|u\|_{D(A^{3/2})}\\
 \label{conto23}
   & \le c \big(1+\|u\|_V^2\big)
   \|u\|_V^{1/2} \|u\|^{3/2}_{D(A^{3/2})}
     \le \delta \|u\|^2_{D(A^{3/2})}
    + c_\delta\|u\|^2_{V}
    + c_\delta\|u\|^{10}_{V},
\end{align}
%
%
%
%
%
%
%
%
%
%
%
%
whence, 
choosing $\delta$ small enough,
we deduce that
\beeq{conto24}
  \|U\|_2
   \le c\big(1+\|V\|_0^5\big), \qquad\textrm{ a.e. in }\,(0,T).
\end{equation}
Let us now recall the Br\'ezis-Gallouet interpolation
inequality \cite[Lemma~2]{BG}, holding
for all $R>0$ and $z\in D(A)$:
\beeq{interpo}
  \|z\|_{L^\infty(\Omega)}
   \le c\OO\|z\|_V\log^{1/2}(1+R)
    +c\OO\|z\|_{D(A)}(1+R)^{-1}.
\end{equation}
%
To apply the inequality, let then $\lambda_0>0$ depending on
$\Omega$ be such that
\beeq{autovalore}
  \|z\|_{D(A)}\ge 2\lambda_0\|z\|_V
   \quad\perogni z\in D(A).
\end{equation}
Then, in case $\lambda_0\|z\|_V\ge 1$, take
$1+R=\|z\|_{D(A)}/\lambda_0\|z\|_V$ in \eqref{interpo}, getting
\bealo
  \|z\|_{L^\infty(\Omega)}
   & \le c\OO\|z\|_V\log^{1/2}\Big(\frac{\|z\|_{D(A)}}{\lambda_0\|z\|_V}\Big)
    +c\OO\|z\|_V
  \le c\OO\|z\|_V\log^{1/2}\|z\|_{D(A)}
    +c\OO\|z\|_V\\
 \label{interpo2}
  & \le c\OO\|z\|_V\log^{1/2}\big(1+\|z\|_{D(A)}\big)
    +c\OO\|z\|_V.
\end{align}
Otherwise, simply choose $R=\|z\|_{D(A)}$, so that
\bealo
  \|z\|_{L^\infty(\Omega)}
   & \le c\OO\|z\|_V\log^{1/2}\big(1+\|z\|_{D(A)}\big)
    +c\OO\frac{\|z\|_{D(A)}}{1+\|z\|_{D(A)}}\\
 \label{interpo3}
   & \le c\OO\|z\|_V\log^{1/2}\big(1+\|z\|_{D(A)}\big)+c\OO.
\end{align}
%
%
%
%
Then, using \eqref{f3}, interpolation,
and either \eqref{interpo2} or \eqref{interpo3},
the \rhs\ of \eqref{conto21new} can be controlled as
follows
\bealo
  \frac12\big(f''(u)v^2,v\big)
   & \le c\big(1+\|u\|_{L^\infty(\Omega)}\big)\|v\|_{L^3(\Omega)}^3\\
 \no
  & \le c\Big(1+\|u\|_V+\|u\|_V
        \log^{1/2}\big(1+\|u\|_{D(A)}\big)\Big)
    \|v\|_{V'}\|v\|_{V}^2\\
 \label{conto26}
  & \le c\Big(1+\cc+\cc
        \log^{1/2}\big(1+\|u\|_{D(A)}\big)\Big)
    \|v\|_{V'}\|v\|_{V}^2.
\end{align}
Thus, estimating the norm of $u$ in
$D(A)$ with the help of \eqref{conto24}
and recalling \eqref{conto25}, inequality \eqref{conto21new}
becomes
\bealo
  \ddt\calF
   +2\kappa\calF
  & \le \cc+c\Big(1+\cc+\cc\log^{1/2}\big(1+\calF\big)\Big)
    \|v\|_{V'}\calF\\
 \label{conto21ter}
  & \le \cc+\kappa\calF
   +c\|v\|_{V'}^2\calF
   \Big(1+\cc+\cc\log\big(1+\calF\big)\Big).
\end{align}
Therefore, possibly replacing $\calF$ with $\calF+c$ for a
suitable $c$ and substituting the expression for $\cc$
from \eqref{st12} (recall that $\cc$ is the quantity
on the \rhs\ and note that if the $\kappa$'s in
\eqref{st12} and \eqref{conto21ter} do not coincide
we can take the smaller),
relation above (for the new $\calF$)
takes the form
\beeq{conto21quater}
  \ddt\calF
   +\kappa\calF
   \le (1+\calF\log\calF)
   \big( Q(\|U_0\|_0)e^{-\kappa t} + Q(\|g\|_{V'}) \big),
\end{equation}
whence the standard theory of ODEs, together with \eqref{conto24},
implies that there exists a computable function
$\calQ:(\RR^+)^3\to \RR^+$, monotone increasing in each of its arguments,
such that
\beeq{stuututt}
  \|u\|_{L^\infty(0,t;D(A^{3/2}))}
   +\|u_t\|_{L^\infty(0,t;V)}
   +\|u_{tt}\|_{L^\infty(0,t;V')}\le
    \calQ\big(\|U_0\|_2,\|g\|_{V'},t\big).
\end{equation}
By standard tools, the above estimate permits to remove the
Faedo-Galerkin approximation and to pass to the limit.
In particular, the regularities~\eqref{regouforte}
are obtained.
This proves the existence part of Theorem~\ref{teoesiforte}.

To prove uniqueness, let us consider a couple of solutions
$u_1$, $u_2$ to (P) in the above regularity setting
and for the same initial data, write
\eqref{CH} for $u_1$ and $u_2$, take the difference, and
test it by $A^{-1} u_t$, where $u:=u_1-u_2$. Using
\eqref{f3} and the regularity \eqref{regouforte},
a straightforward computation gives
\bealo
  \io \big(f(u_1)-f(u_2)\big)u_t
   & \le \frac12\|u_t\|_{V'}^2
     +\frac12\io\big|f'(u_1)\nabla u_1-f'(u_2)\nabla u_2\big|^2\\
  \label{contouniq}
  & \le \frac12\|u_t\|_{V'}^2
     +C\|u\|_V^2,
\end{align}
where $C$ depends, of course, on the norms of $u_1$, $u_2$
specified in \eqref{regouforte}. Thus, Gronwall's lemma permits
to conclude that $u_1\equiv u_2$ as desired, which completes
the proof of Theorem~\ref{teoesiforte}.\dimbox
%
%
%
%

      \section{Asymptotic behavior of quasi-strong solutions}
       \label{secV2}


We associate with Problem~(P) the semiflow $\calS$ acting
on $\calV_2$ and generated by the
quasi-strong solutions provided  by Theorem~\ref{teoesiforte}.
We will also indicate by $S(t)$,
$t\ge 0$, the semigroup operator defined by $\calS$.
Let us now prove some important
properties of $\calS$ and $S(t)$.
%
%
\bete\label{teosemi}
 Let the assumptions of\/ {\rm Theorem~\ref{teoesiforte}}
 hold. Then, the semiflow~$\calS$ is\/ {\rm uniformly
 dissipative}. Namely, there exists a constant $R_0$
 independent of the initial data such that, for all
 bounded $B\subset \calV_2$, there exists $T_B\ge 0$
 such that $\|S(t)b\|_2\le R_0$, for all
 $b\in B$ and $t\ge T_B$.
 Moreover, any $u\in \calS$ satisfies the additional
 time continuity property
 \beeq{regouforte2}
   u\in C^2([0,T];V')\cap C^1([0,T];V)\cap C^0([0,T];D(A^{3/2})).
 \end{equation}
 Finally, given a sequence of initial data
 $\{(u\znn,u\unn)\}\subset \calV_2$ suitably tending to
 some $(u_0,u_1)\in \calV_2$,
 and denoting by $u_n,u$ the solutions
 emanating from $(u\znn,u\unn),(u_0,u_1)$, respectively,
 we have that
 \beal{conticompaw}
   (u\znn,u\unn)\to (u_0,u_1) \text{~~weakly in }\/\calV_2
    & ~\Rightarrow~ (u_n,u_{n,t}) \to (u,u_{t})
    \text{~~weakly star in }\,L^\infty(0,T;\calV_2),\\
 \label{conticompa}
   (u\znn,u\unn)\to (u_0,u_1) \text{~~strongly in }\/\calV_2
    & ~\Rightarrow~ (u_n,u_{n,t}) \to (u,u_{t})
    \text{~~strongly in }\,C^0([0,T];\calV_2),
 \end{align}
 for any fixed $T\ge 0$.
\ente
\noindent%
\begin{proof}
Let us prove the existence of a bounded absorbing set first. Coming back
to \eqref{conto21ter}, we can now rewrite it in the form
\beeq{conto21quinquies}
  \ddt\calF
   +\kappa\calF
   \le \big(1+\|v\|_{V'}^2\calF\log\calF\big)
   \big( Q(\|U_0\|_0)e^{-\kappa t} + Q(\|g\|_{V'}) \big).
\end{equation}
Next, recalling \eqref{st12} and \eqref{dissin},
we can take $T_1$ so large,
only depending on $\|U_0\|_0$, that
\beal{partenzagronw1}
  & Q(\|U_0\|_0)e^{-\kappa t} + Q(\|g\|_{V'})
   \le c_1, \qquad\perogni t\ge T_1\\
 \label{partenzagronw2}
  & \int_t^\infty\|v(s)\|_{V'}^2\,\dis\le c_2,
    \qquad\perogni t\ge T_1,
\end{align}
where $c_1$, $c_2$ depend on $g$ but
do not depend on $U_0$.
Thus, setting $y:=\log \calF\ge 0$,
\eqref{conto21quinquies} can be rewritten,
for $t\ge T_1$, as
\beeq{conto21sexies}
  y'+\kappa
   \le c_1 e^{-y} + c_1\|v\|_{V'}^2y.
\end{equation}
Moreover, by \eqref{stuututt}, and possibly modifying the
expression of $\calQ$, we have
\beeq{FT1}
  y(T_1) \le \log\big(\calQ(\|U_0\|_2,
    \|g\|_{V'},T_1)\big)
   =:\eta,
\end{equation}
where the value of $\eta$ depends only on
$\|U_0\|_2$ since so does $T_1$.
Assume now that $c_1>\kappa$
(if not, we can suitably modify its value).
Setting $\zeta:=\log(2c_1/\kappa)>0$,
we distinguish
\beeq{whether}
  \text{whether }\,y(T_1)>\zeta
   \quext{or }\,y(T_1)\le \zeta.
\end{equation}
If the first condition holds, in
a right neighborhood of $T_1$ it is
\beeq{gronw2}
  \ddt\Big(y+\frac\kappa2 t\Big)
   =y'+\frac\kappa2
   \le c_1\|v\|_{V'}^2y
   \le c_1\|v\|_{V'}^2\Big(y+\frac\kappa2 t\Big),
\end{equation}
so that, recalling \eqref{partenzagronw2}
and \eqref{FT1}, we obtain
%
\beeq{solvlin1}
   y(t) 
    \le \Big(\eta
     +\frac\kappa2 T_1\Big) e^{c_1c_2}-\frac{\kappa}2 t,
\end{equation}
which implies that for some time $\tau_1>T_1$,
still depending only on $\|U_0\|_2$,
$y(\tau_1)=\zeta$.

Thus, we have essentially
reduced us to the case when $y(T_1)\le \zeta$,
which we now treat. Let us then set
$\zeta_*:=\log(c_1/\kappa)>0$ and define
$y_*:=y\vee \zeta_*$. It is clear that
$y_*$ satisfies $y_*(T_1)\le \zeta$ and,
for almost all~$t\ge T_1$,
\beeq{disugy*}
  y_*' \le c_1 y_*\|v\|_{V'}^2,
\end{equation}
so that, solving \eqref{disugy*},
noting that $y\le y_*$,
and using \eqref{partenzagronw2},
we have
\beeq{disugy*bis}
  y(t)\le y_*(t)\le \zeta\exp(c_1c_2),
   \quext{respectively }\,\perogni t\ge T_1
   \quext{or }\,\perogni t\ge \tau_1.
\end{equation}
Thus, we can conclude the proof of
the ``dissipative'' part of Theorem~\ref{teosemi}
by taking $T_B:=\max\{T_1,\tau_1\}$ (where the choice
of $T_1$ or $\tau_1$ corresponds now to the radius
of the chosen bounded set $B\subset\calV_2$). Note that
the explicit value of the radius $R_0$
of the absorbing ball, i.e., the \rhs\ of \eqref{disugy*bis},
could be explicitly computed by referring
to \eqref{conto25} and \eqref{conto24}.

\vspace{2mm}

\noindent%
Next, let us show \eqref{regouforte2}.
With this aim, let us first rewrite \eqref{CH} in the form
\beeq{CH2}
  u\dtt+u_t+A^2u=f'(u)\Delta u+f''(u)|\nabla u|^2+g=:G+g,
\end{equation}
which makes sense thanks to \eqref{regouforte} and
\eqref{f1}-\eqref{f3}. Actually, one can also easily
prove that
\beeq{regoG}
  G\in L^2(0,T;V)
\end{equation}
(much more is true, in fact, but the above
is sufficient for what follows). Then,
given some $T>0$, we take sequences such that
\beal{appro1}
  \big\{G_n\big\}\subset C^0([0,T];D(A)),
   & \qquad G_n\to G
   \quext{strongly in }\,L^2(0,T;V),\\
 \label{appro2}
  \big\{(u\zzn,u_{1,n})\big\}\subset \calV_3,
   & \qquad (u\zzn,u_{1,n})\to(u_0,u_1)
   \quext{strongly in }\,\calV_2,\\
 \label{appro3}
  \big\{g_n\big\}\subset V,
   & \qquad g_n\to g
   \quext{strongly in }\,V',
\end{align}
and, for all $n\in\NN$, we consider the
solution $U_n=(u_n,u_{n,t})$ to
\beeq{CHn}
  u_{n,tt}+u\nnt+A^2u_n=G_n+g_n,
\end{equation}
coupled with the new initial datum
$U_{0,n}=(u\zzn,u_{1,n})$. By the linear theory,
this satisfies
\beeq{regounew3}
  u_n\in C^2([0,T];H)\cap C^1([0,T];D(A))\cap C^0([0,T];D(A^2)),
\end{equation}
so it is suitable for the a priori estimates
we need. Writing \eqref{CHn} for the couple of
indexes $n,m$, taking the difference, temporarily
setting $u:=u_n-u_m$, testing by $A u_t$,
and integrating over $(0,t)$ for $t\le T$,
we readily get
\bealo
  & \|U(t)\|_2^2
   -2\duavg{g_n-g_m,Au(t)}
   +\itt\|\nabla u_t\|^2\\
 \label{conto31}
  & \mbox{}~~~~~
  \le\big\|U_{0,n}-U_{0,m}\big\|_2^2
   -2\duavg{g_n-g_m,A(u_{0,n}-u_{0,m})}
   +\itt\|\nabla (G_n-G_m)\|^2.
\end{align}
Thus, taking the supremum with respect to~$t\in [0,T]$
we deduce that $\{u_n\}$ is a Cauchy sequence
with respect to the latter two norms in \eqref{regouforte2}.
Since the convergence of $u_{tt}$ can be proved
by a comparison of terms in \eqref{CHn},
this entails that $u$ fulfills \eqref{regouforte2}.

\vspace{2mm}

To conclude, let us examine the continuity properties
of $S(t)$. First of all, \eqref{conticompaw}
can be shown by using uniform boundedness in
the norms indicated there, weak compactness,
and lower semicontinuity of norms with respect to~weak star convergence.
Note that we do not need to extract subsequences
since we have uniqueness of the limit.

The proof of \eqref{conticompa} is equally simple but
more technical. For this reason, we proceed by
deriving formal estimates and just give the highlights
of the procedure which could be used to make them rigorous.
Thus, let us write \eqref{CH} for $u_n$, $u_m$
($u_n$, $u_m$ being now as in the statement),
take the difference, and differentiate the resulting equation
with respect to time. Testing by $A^{-1}u_{tt}$,
where $u:=u_n-u_m$, one then formally infers
\bealo
  & \ddt\|U_t\|_0^2
   +\|u_{tt}\|_{V'}^2
   \le \big\|f'(u_n)u_{n,t}-f'(u_m)u_{m,t}\big\|_{V}^2\\
 \label{conto41}
  & \mbox{}~~~~~
  \le 2\big\|f''(u_n)\nabla u_n u_{n,t}-f''(u_m) \nabla u_m
       u_{m,t}\big\|^2
   +2\big\|f'(u_n)\nabla u_{n,t}-f'(u_m) \nabla u_{m,t}\big\|^2,
\end{align}
whence performing standard calculations (i.e., adding and
subtracting some terms, exploiting uniform boundedness
in the norms specified in \eqref{regouforte},
and using the properties of $f$ as well as suitable
Sobolev's embeddings),
one can transform \eqref{conto41} into
\beeq{conto42}
  \ddt\|U_t\|_0^2
   +\|u_{tt}\|_{V'}^2
   \le C_0\big(
   \|u_t\|_V^2+\|u\|_{D(A)}^2\big),
\end{equation}
where $C_0$ depends on the norms in
\eqref{regouforte} of $u_n$ and $u_m$. 

Next, let us write again \eqref{CH} for $u_n$, $u_m$, take
the difference, and test it by $u$.
This yields
\beeq{conto43}
  \ddt\|u\|^2
   +\|u\|_{D(A)}^2
  \le c\big(\|u_{tt}\|_{V'}^2
   +\|f(u_n)-f(u_m)\|^2\big)
  \le C\big(\|u_{tt}\|_{V'}^2
   +\|u\|^2\big),
\end{equation}
whence, multiplying \eqref{conto43} by $2C_0$,
and adding the result to \eqref{conto42},
an application of Gronwall's lemma leads to the
strong convergences of $u_{n,tt}$ to $u_{tt}$
in $C^0([0,T];V')$ and of $u_{n,t}$ to $u_{t}$
in $C^0([0,T];V)$. Note, indeed, that \eqref{conticompa}
and a comparison in \eqref{CH} guarantee that
$u_{n,tt}(0)\to u_{tt}(0)$ strongly in $V'$.
Finally, the strong convergence
of $u_{n}$ to $u$ in $C^0([0,T];D(A^{3/2}))$
can be proved by a further comparison of terms
in the (difference) of \eqref{CH}.

Of course, the above procedure is not fully rigorous
since the test function $A^{-1}u_{tt}$ is not
admissible for the time derivative of \eqref{CH}.
To overcome this problem, one could argue as in the
proof of \eqref{regouforte2}. Namely, setting
\beeq{Gnm}
  G_{n,m}:=f(u_n)-f(u_m),
\end{equation}
for fixed $m$ and $n$, one notes that
\beeq{regoGnm}
  G_{n,m}\in C^1([0,T];V)\cap C^0([0,T];D(A))
\end{equation}
and can approximate $G_{n,m}$ by
a sequence $\{G_{n,m}^k\}\subset C^1([0,T];D(A^{3/2}))$ such that
\beeq{Gnmk}
   G_{n,m}^k\xrightarrow{k} G_{n,m} \quext{strongly in }\,
   C^1([0,T];V)\cap C^0([0,T];D(A)).
\end{equation}
Then, suitably approximating also the initial and source
data (cf.~\eqref{appro2}--\eqref{appro3}), and noting that
the $k$-solution $u_k$ is sufficiently regular,
one can perform the estimates described above working
on $u_k$ and then take the limit with respect to~$k$.
The details are left to the reader. The proof
is complete.
\end{proof}

\vspace{2mm}

\noindent%
The next theorem,
whose proof relies on the so-called ``energy
method'' (cf.~\cite[Sec.~2]{Ba2} for the theoretical
background and a comparison with the second order case, see also~\cite{MRW}),
states the asymptotic compactness in $\calV_2$ of~$\calS$.
\bete\label{teoasy}
 Let the assumptions of\/ {\rm Theorem~\ref{teoesiforte}}
 hold. Then, the semiflow $\calS$ associated
 to~{\rm (P)} is\/ {\rm asymptotically compact}.
 Namely, for any $\calV_2$-bounded
 sequence $\{(u\znn,u\unn)\}$ of initial data
 and any positively diverging sequence
 $\{t_n\}$ of times, there
 exists $(\chi,\chi_1)\in \calV_2$ such that a subsequence
 of $\{(u_n(t_n),u_{n,t}(t_n))\}$ tends to $(\chi,\chi_1)$\/
 {\rm strongly} in~$\calV_2$ ($u_n$ is here the solution
 having $(u\znn,u\unn)$ as initial datum).
\ente
\noindent%
\begin{proof}
Let us first notice that, as a consequence of
the first part of the previous proof (cf., in
particular, \eqref{conto31}), the solution $u$ to~(P)
satisfies for all $s,t$ the {\sl equality}
\beeq{quasienergy}
  \frac12\|U(t)\|_2^2
   -\duav{g,Au(t)}
   +\int_s^t\|\nabla u_t\|^2
  =\frac12\|U(s)\|_2^2
   -\duav{g,Au(s)}
   +\int_s^t\big\langle Au_t,G\big\rangle,
\end{equation}
where $G$ has been defined in \eqref{CH2}
and $U=(u,u_t)$.
Next, as we substitute $G$ with its expression,
we claim that, for all $s,t\in [0,T]$,
\beeq{ipepan}
  \int_s^t\big\langle A u_t,-f'(u)\Delta u\big\rangle
   =\frac12\io f'(u(t))|\Delta u(t)|^2
     -\frac12\io f'(u(s))|\Delta u(s)|^2
   -\frac12\int_s^t\io f''(u)u_t|\Delta u|^2.
\end{equation}
To prove this, let us proceed once more by regularization.
Actually, \eqref{ipepan} surely holds for a
more regular $u_n$. Assuming that
\beeq{counu}
  u_n\to u \quext{strongly in }\,C^1([0,T];V)
   \cap C^0([0,T];D(A^{3/2})),
\end{equation}
we can write \eqref{ipepan} for $u_n$ and take the limit. This is
straightforward as far as the \rhs\ is concerned.
On the other hand, the integrand on the \lhs\
can be rewritten as
\beeq{ipepa2}
  -\big(\nabla u\nnt,\nabla (f'(u_n)\Delta u_n)\big),
\end{equation}
and it is easy to prove,
using \eqref{counu}, that both
terms in the scalar product converge to the expected
limits in the strong topology of $\CZH$.
Since
\beeq{ipepa2bis}
  -\big(\nabla u_t,\nabla (f'(u)\Delta u)\big)
   =\big\langle A u_t,-f'(u)\Delta u\big\rangle
\end{equation}
by definition of $A$, this concludes
the proof of \eqref{ipepan}.

Consequently, $u$ turns out to satisfy
the following higher order energy equality
\bealo
  & \frac12\|U(t)\|_2^2
   -\duav{g,Au(t)}
   +\int_s^t\|\nabla u_t\|^2
   +\frac12 \io f'(u(t))|\Delta u(t)|^2\\
 \no
  & \mbox{}~~~~~
  =\frac12\|U(s)\|_2^2
   -\duav{g,Au(s)}
   +\frac12 \io f'(u(s))|\Delta u(s)|^2\\
 \label{energy}
  & \mbox{}~~~~~~~~~~~~~~~
  +\frac12\int_s^t\io f''(u)u_t|\Delta u|^2
  +\int_s^t \big\langle{Au_t,f''(u)|\nabla u|^2}\big\rangle.
\end{align}
This is the starting point to implement the
so-called ``energy method'' introduced by J.M.~Ball
(cf.~\cite[Sec.~4]{Ba2}, see also \cite{MRW})
to prove asymptotic compactness,
which is our next task.

To start with, let us define the functional
\beeq{deficalG0}
  \calG_0(t)
   :=\frac12\|U(t)\|_2^2
   -\duav{g,Au(t)}
   +\frac12 \io f'(u(t))|\Delta u(t)|^2.
\end{equation}
Actually, at least if no danger of confusion occurs,
we shall write indifferently $\calG_0$,
$\calG_0(t)$ or $\calG_0(u(t))$ in the sequel,
with some abuse of language since in fact
$\calG_0$ depends both on $u$ and on $u_t$.
We shall use the same convention also for the
other functionals defined below.

Then, writing \eqref{energy} for $t=s+h$, dividing by $h$,
and letting $h\to 0$,
it is immediate to deduce that $\calG_0$ is absolutely
continuous on $[0,T]$ and there holds
\beeq{energydiff}
  \ddt\calG_0
   +\|\nabla u_t\|^2
   =\frac12\io f''(u)u_t|\Delta u|^2
   +\big\langle{Au_t,f''(u)|\nabla u|^2}\big\rangle,\quad\textrm{ a.e. in }\,(0,T).
\end{equation}
Next, let us test \eqref{CH} by $Au$. The same procedure
used before permits to justify the validity, a.e.~in $(0,T)$,
of the equality
\beeq{energydiff2}
  \ddt\Big(
   (u_t,Au)
    +\frac12\|\nabla u\|^2\Big)
    -\|\nabla u_t\|^2
   +\|u\|_{D(A^{3/2})}^2
   +\io f'(u)|\Delta u|^2
   -\duav{g,Au}
  =-\io f''(u)|\nabla u|^2\Delta u.
\end{equation}
Then, let us multiply \eqref{energydiff2} by $1/2$ and
sum the result to \eqref{energydiff}. We get
\beeq{energydiff2.2}
  \ddt\Big(\calG_0+\frac12(u_t,Au)
    +\frac14\|\nabla u\|^2\Big)
    +\frac12\|U\|_2^2
    +\frac12\io f'(u)|\Delta u|^2
   -\frac12\duav{g,Au}
  =\calH_0,
\end{equation}
where we have set
\beeq{deficalH0}
  \calH_0:=\frac12\io f''(u)u_t|\Delta u|^2
   +\big\langle{Au_t,f''(u)|\nabla u|^2}\big\rangle
   -\frac12\io f''(u)|\nabla u|^2\Delta u.
\end{equation}
Consequently, adding some terms to both hands sides of
\eqref{energydiff2.2} we obtain the equality
\beeq{energydiff3}
  \ddt\calG+\calG=\calH, \quext{a.e.~in }\,(0,T),
\end{equation}
where we have set
\beal{deficalG}
  \calG & := \calG_0
    +\frac12(u_t,Au)
   +\frac14\|\nabla u\|^2,\\
 \label{deficalH}
  \calH & := \calH_0
   -\frac12\duav{g,Au}
   +\frac12(u_t,Au)
   +\frac14\|\nabla u\|^2.
\end{align}
Thus, from \eqref{energydiff3} and for any $\tau,M\ge 0$,
we obtain
\beeq{energyint}
  \calG(\tau+M)=\calG(\tau)e^{-M}
   +\int_\tau^{\tau+M}e^{s-\tau-M}\calH(s)\,\dis.
\end{equation}
At this point, recalling the notation
in the statement, let
us set $v_n(t):=u_n(t_n+t-M-\tau)$ (so that, in particular,
$v_n(\tau)=u_n(t_n-M)$ and $v_n(\tau+M)=u_n(t_n)$).
Since $(u_{0,n},u_{1,n})$ is bounded in $\calV_2$,
by uniform dissipativity it follows that
$u_n$ is uniformly bounded in the norms
\eqref{regouforte} (in a way which does not depend
on $T$) by some constant $C$.
The same, of course, holds also for $v_n$
and for the values of the functional $\calG$.
Thus, at least a (nonrelabelled)
subsequence of $v_n$ tends to a solution
$v$ weakly star in the norm specified in
\eqref{regouforte} and for all $T>0$.
More precisely, we have
\beeq{newdebo}
  v_n\to v \quext{in }\,C^1_w([0,T];V)
      \cap C^0_w([0,T];D(A^{3/2})),
\end{equation}
so that, in particular, there exist
the limits
\beal{defichi}
  & (\chi,\chi_1):=\lim_{n\to\infty} \big(v_n(\tau+M),v_{n,t}(\tau+M)\big)
     =\lim_{n\to\infty} \big(u_n(t_n),u_{n,t}(t_n)\big),\\
 \label{defichi2}
 &  (\chi_{-M},\chi_{1,-M}):=\lim_{n\to\infty} \big(v_n(\tau),v_{n,t}(\tau)\big)
     =\lim_{n\to\infty} \big(u_n(t_n-M),u_{n,t}(t_n-M)\big),
\end{align}
which have to be intended, at least in the meanwhile,
in the {\sl weak topology}\/
of $\calV_2$. Moreover, it is
$(\chi,\chi_1)=(v(\tau+M),v_t(\tau+M))$ and
$(\chi_{-M},\chi_{1,-M})=(v(\tau),v_t(\tau))$.
At this point, writing \eqref{energyint}
for $v_n$ gives
\beeq{energyintvn}
  \calG(v_n(\tau+M))-\calG(v_n(\tau))e^{-M}
   =\calG(u_n(t_n))-\calG(u_n(t_n-M))e^{-M}
   =\int_\tau^{\tau+M}e^{s-\tau-M}\calH(v_n(s))\,\dis.
\end{equation}
It is now a standard procedure
to check that, as far as $\tau,M$ are fixed,
\eqref{newdebo} and suitable Sobolev's embeddings
give
\beeq{conveHn}
  \int_\tau^{\tau+M}e^{s-\tau-M}\calH(v_n(s))\,\dis
   \to\int_\tau^{\tau+M}e^{s-\tau-M}\calH(v(s))\,\dis.
\end{equation}
Thus, taking the supremum limit of \eqref{energyintvn}
as (a proper subsequence of) $n$ tends to $\infty$,
we get
\bealo
  \limsup_{n\nearrow\infty}
   \calG(u_n(t_n))
  & \le C e^{-M}
   +\limsup_{n\nearrow\infty}
       \int_\tau^{\tau+M}e^{s-\tau-M}\calH(v_n(s))\,\dis\\
 \no
  & = C e^{-M}
    +\int_\tau^{\tau+M}e^{s-\tau-M}\calH(v(s))\,\dis\\
 \no
  & = C e^{-M}
   + \calG(v(\tau+M))-\calG(v(\tau))e^{-M}\\
 \label{contoball}
  & \le C e^{-M}
   + \calG(\chi),
\end{align}
where in deducing the third equality
we used that \eqref{energyint}
is satisfied also by the limit solution $v$.
Since the above holds for all $M>0$ and
with $C$ independent of $M$, letting $M\nearrow\infty$
and using the immediate fact
that $\calG$ is sequentially weakly lower semicontinuous
in $\calV_2$, we eventually obtain that
\beeq{coG}
  \calG(u_n(t_n))\to \calG(\chi).
\end{equation}
Thus, since it is clear that
just the weak convergence
\eqref{defichi} entails
\bealo
 &  -\duavg{g,Au_n(t_n)}
   +\frac12\io f'(u_n(t_n))|\Delta u_n(t_n)|^2
   +\frac12\big(u_{n,t}(t_n),Au_n(t_n)\big)
   +\frac14\big\|\nabla u_n(t_n)\big\|^2\\
\label{lowerorder}
  & \mbox{}~~~~~ \to
   -\duavg{g,A\chi}
   +\frac12\io f'(\chi)|\Delta \chi|^2
   +\frac12(\chi_1,A\chi)
   +\frac14\|\nabla \chi\|^2,
\end{align}
comparing \eqref{coG} with \eqref{lowerorder} we obtain
\beeq{cofinal}
  \big\|(u_n(t_n),u_{n,t}(t_n))\big\|_2\to
   \|(\chi,\chi_1)\|_2.
\end{equation}
This relation, together with the weak convergence
\eqref{defichi}, gives
the desired strong convergence and concludes the proof
of asymptotic compactness and of
Theorem~\ref{teoasy} as well.
\end{proof}

On account of Theorems \ref{teosemi} and \ref{teoasy}, we thus deduce (see, e.g., \cite[Thm.~3.3]{Ba1})
\bete\label{teoattrafirst}
 Let the assumptions of\/ {\rm Theorem~\ref{teoesiforte}} hold.
 Then, the semiflow~$\calS$ possesses the global attractor~$\calA_2$.
\ente


      \section{Smoothness of the global attractor $\calA_2$}
       \label{secV2b}


We prove here a regularity property of the attractor
constructed in the previous section. This fact will
be used for constructing an
exponential attractor in Section~\ref{secexpo}.
A straightforward consequence of the results of this section is
the existence of the global attractor which was already obtained
in Theorem~\ref{teoasy} by a different technique.
We decided to keep the latter proof because it is simple and (hopefully)
interesting in itself. On the contrary, the proof of Theorem~\ref{teoattrastrong},
which relies on a new decomposition method partly
related to that in \cite{PZ}, involves a number of technical complications.
We also point out that, in Section~\ref{secV0}, we will appeal to
the same technique used for proving Theorem~\ref{teoasy} to establish the existence of the
global attractor for weak solution (see Theorem~\ref{teoasydebo}).

%
\bete\label{teoattrastrong}
 Let the assumptions of\/ {\rm Theorem~\ref{teoesiforte}}
 hold. Additionally, let
 \beeq{regogplus}
   g\in H.
 \end{equation}
 Then, the global attractor~$\calA_2$ for
 the semiflow~$\calS$ is bounded
 in $\calV_3$.
\ente
\beos\label{attr2stru}
On account of well-known results, we can infer that $\calA_2$ consists
of those points of $\calV_3$ from which bounded complete trajectories
originate. These are strong solutions to\/ {\rm Problem~(P)}.
\eddos
\noindent
{\bf Proof of Theorem~\ref{teoattrastrong}.}
The proof is divided into several steps which are
presented as separate lemmas. We start with a simple
property whose proof is more or less straightforward.
\bele\label{lemmaR2}
 Let $S:[0,+\infty)^2\to [0,\infty)$, $S=S(t,R)$,
 be a continuous function such that\\[2mm]
 (i)~~There exists $R_0\in [0,\infty)$
 such that for all $R\in [0,\infty)$
 there exists $T_R\in [0,\infty)$
 such that $S(t,R)\le R_0$
 for all $t\ge T_R$;\\[1mm]
 (ii)~~$S(0,R)=R$ for all $R\in [0,\infty)$;\\[1mm]
 (iii)~~$R\mapsto S(t,R)$ is increasingly monotone for
 all $t\in [0,\infty)$.\\[2mm]
 Then, there exists an increasingly monotone
 function $Q:[0,\infty)\to[0,\infty)$,
 $Q=Q(R)$, such that $S(t,R)\le Q(R) e^{-t}+R_0$ for
 all $(t,R)\in [0,\infty)^2$.\dimbox
\enle
Actually, noting as $S(t,R)$ the $\calV_2$-radius
of $S(t)B(0,R)$ ($B(0,R)$ being the $R$-ball in $\calV_2$),
it is clear that $S$ verifies the
properties (i)--(iii). Thus, using Lemma~\ref{lemmaR2},
the dissipativity property of Theorem~\ref{teosemi} is
rewritten as
\beeq{dissifo}
  \|S(t)U_0\|_2\le Q\big(\|U_0\|_2\big)e^{-t}+R_0,
   \quext{where }\,R_0=Q(\|g\|_{V'}).
\end{equation}
Next lemma states that $\calV_2$-solutions to~(P)
satisfy a dissipation property similar
to \eqref{dissin}, but in a stronger norm.
\bele\label{lemmaZelik1}
 Let the assumptions of\/ {\rm Theorem~\ref{teoattrastrong}} hold.
 Then, for any $\calV_2$-solution $U=(u,u_t)$ to~{\rm (P)}
 there holds
 \beeq{relZelik1}
   \int_0^\infty\|u_t(s)\|^2_{V}\,\dis
    +\sup_{s\in[0,\infty)}\|u_{tt}(s)\|^2_{V'}
   \le Q\big(\|U_0\|_2\big)<\infty.
 \end{equation}
\enle
\begin{proof}
 Let us differentiate \eqref{CH} with respect to time and set
 $\theta:=u_t$ and $\Theta:=(\theta,\theta_t)$.
 We get, for $L>0$ to be chosen later,
 \beeq{CHt}
   A^{-1}(\theta_{tt}+\theta_t)
     +A\theta+LA^{-1}\theta+f'(u)\theta=h:=LA^{-1}u_t.
 \end{equation}
 We now test \eqref{CHt} by $2\theta_t$.
 Simple computations lead to
 \bealo
   & \ddt\big[ \|\Theta\|_0^2
    + L \|\theta\|_{V'}^2
    + \big( f'(u) \theta, \theta \big)\big]
      + \|\theta_t\|_{V'}^2 \\
  \label{Zelik11}
   & \mbox{}~~~~~
    \le \|h\|_{V}^2 + \big( f''(u) \theta^2, \theta \big)
    \le \|h\|_{V}^2 + \frac14 \|\theta\|_V^2 + C \|\theta\|_{V'}^2,
  \end{align}
 where we have used the duality pairing
 $V'$-$V$, Young's inequality, and, in the last
 passage, the uniform $\calV_2$-boundedness of
 $U$ and interpolation. Note that, here
 and in the rest of this Section, the constants
 $C$ are allowed to depend (actually, at
 most in a polynomial way since $f$ grows
 polynomially) on the $\calV_2$-norm of the solution,
 which is uniformly bounded in time by \eqref{dissifo}.
 In fact, $C$ is a quantity having the same expression
 as the \rhs\ of \eqref{dissifo}, but
 we do {\sl not}\/ allow $C$ to
 depend on $L$.

 Next, we test \eqref{CHt} by $\theta/2$, inferring
 \bealo
   & \ddt\Big[ \frac12 \big(\theta_t,A^{-1}\theta\big)
    + \frac14\|\theta\|_{V'}^2 \Big]
    -\frac12 \|\theta_t\|_{V'}^2
    +\frac12 \|\theta\|_{V}^2
    +\frac{L}2 \|\theta\|_{V'}^2
    +\frac12 \big( f'(u) \theta, \theta \big)\\
 \label{Zelik12}
  & \mbox{}~~~~~
   = \frac12 \big(h,\theta\big)
   \le \frac18 \|\theta\|^2_V
   + c \|h\|^2_{V}
 \end{align}
 Summing \eqref{Zelik11} and \eqref{Zelik12},
 we then get
 \bealo
   & \ddt\Big[ \|\Theta\|_0^2
    + \Big(L+\frac14\Big) \|\theta\|_{V'}^2
    + \big( f'(u) \theta, \theta \big)
    +\frac12 \big(\theta_t,A^{-1}\theta\big) \Big] \\
  \label{Zelik13}
   & \mbox{}~~~~~~~~~~
    + \frac12\|\theta_t\|_{V'}^2
    +\frac18 \|\theta\|_{V}^2
    +\frac{L-2C}2 \|\theta\|_{V'}^2
    +\frac12 \big( f'(u) \theta, \theta \big)
   \le c\|h\|_{V}^2,
 \end{align}
 still for $C$ independent of $L$. Thus,
 noting by $\calY$ the quantity in square brackets,
 we notice that we can choose $L$ so large (depending
 on $C$ and $\lambda$ in \eqref{f3}) to get
 \beeq{Zelik14}
   \ddt \calY + \kappa \calY \le c\|h\|_{V}^2,
 \end{equation}
 for some $\kappa>0$. Hence, recalling that
 $h=LA^{-1}u_t$ and using \eqref{dissin}
 and the fact that $\calY(0)=Q(\|U_0\|_2)$,
 \eqref{relZelik1} follows immediately.
\end{proof}
We are now ready to decompose the solution $u$
to~(P) as the sum of a ``compact'' part
\beeq{CHcomp}
 A^{-1}(v_{tt}+v_t)
    +Av+LA^{-1}v+f(v)=LA^{-1}u+g,
 \qquad V|_{t=0}=0,
\end{equation}
where $V:=(v,v_t)$, and a ``decaying'' part
\beeq{CHdecay}
 A^{-1}(w_{tt}+w_t)
    +Aw+LA^{-1}w+f(u)-f(v)=0,
 \qquad W|_{t=0}=U_0:=(u_0,u_1),
\end{equation}
where $W:=(w,w_t)$ and $U_0$ belongs to a bounded absorbing set in $\calV_2$ (cf. \eqref{dissifo}).
Note that the value of $L$ in
\eqref{CHcomp}-\eqref{CHdecay} will possibly differ
from that in \eqref{CHt}.
\bele\label{lemmaZelik2}
 Let the assumptions of\/ {\rm Theorem~\ref{teoattrastrong}} hold.
 Then, $L$ can be chosen so large that
 \beeq{relZelik2}
   \|W(t)\|_0
   \le Q(\|U_0\|_2) e^{-\kappa t}.
 \end{equation}
 %
\enle
\begin{proof}
 We proceed along the lines of the preceding proof.
 First, we test \eqref{CHdecay} by $2 w_t$, so that
 \beeq{Zelik21}
   \ddt\big[ \|W\|_0^2
    + L \|w\|_{V'}^2
    + 2 \calI_1 \big]
    + \|w_t\|_{V'}^2
   \le 2 \calI_2,
 \end{equation}
 where $\calI_1$ and $\calI_2$ collect the terms coming
 from $f$. Namely, we have (recall that $U=V+W$)
 \beeq{stimaI1}
   \calI_1=\big(F(u-w)-F(u)+f(u)w,1\big)
    \ge -\frac\lambda2 \|w\|^2
 \end{equation}
 thanks to \eqref{f2}. To estimate $\calI_2$,
 let us first notice that,
 performing the standard energy estimate
 (cf.~\eqref{conto11}--\eqref{st12}) on
 \eqref{CHcomp} (i.e., testing it by
 $v_t+\delta v$ for small $\delta>0$)
 and using the energy estimate
 \eqref{st12} for $u$ to control the term on
 the \rhs, we derive
 \beeq{dissicompa}
   \|V(t)\|_0^2\le Q_L(\|U_0\|_0) e^{-\kappa t} + C_L,
 \end{equation}
 where both $C_L$ and $Q_L$ depend on $L$ since
 so does the \rhs\ of \eqref{CHcomp}. Comparing
 \eqref{st12} and \eqref{dissicompa}, we also
 get
 \beeq{dissidissi}
   \|W(t)\|_0^2\le Q_L(\|U_0\|_0) e^{-\kappa t} + C_L.
 \end{equation}
 Using now the uniform $\calV_2$-bound on $U$,
 \eqref{f3}, \eqref{dissidissi}, and standard interpolation and
 embeddings, we can estimate
 \bealo
   \calI_2 & = \big(f(u-w)-f(u)+f'(u)w,u_t\big)
     \le c \io \big(1+|u|+|w|\big)
      |w|^2|u_t|\\
  \label{stimaI2}
    & \le C_L \|u_t\|_V \|w\|_V^2
     \le  \frac1{16} \|w\|_V^2 + C_L\|u_t\|_V^2 \|w\|_V^2,
 \end{align}
 where the dependence on $L$ of the constant $C_L$
 comes from \eqref{dissidissi}.

 Next, we test \eqref{CHdecay} by $w/2$, inferring
 \beeq{Zelik22}
   \ddt\Big[ \frac12 \big(w_t,A^{-1}w\big)
    + \frac14\|w\|_{V'}^2 \Big]
    -\frac12 \|w_t\|_{V'}^2
    +\frac12 \|w\|_{V}^2
    +\frac{L}2 \|w\|_{V'}^2
    +\frac12\big(f(u)-f(u-w),w\big)
   =0.
 \end{equation}
 Now, using \eqref{f2} and interpolation,
 it is not difficult to compute
 \bealo
   & \big(f(u)-f(u-w),w\big)
   = \big(f(u)-f(u-w),w\big)
      + \big(F(u-w)-F(u),1\big) - \big(F(u-w)-F(u),1\big)\\
  \label{Zelik23}
  & \mbox{}~~~~~
   = \calI_1
      +\big(F(u)-F(u-w)-f(u-w)w,1\big)
   \ge \calI_1 - \frac\lambda{2} \|w\|^2
   \ge \calI_1 - \frac1{4} \|w\|_V^2 - C \|w\|_{V'}^2,
 \end{align}
 for some (new) $C>0$.
 Thus, summing \eqref{Zelik21} and \eqref{Zelik22}
 we arrive at
 \bealo
   & \ddt\Big[ \|W\|_0^2
    + \Big(L+\frac14\Big) \|w\|_{V'}^2
    + 2 \calI_1
    + \frac12 \big(w_t,A^{-1}w\big) \Big]\\
    \nonumber
   & +\frac12 \|w_t\|_{V'}^2
    +\frac12 \|w\|_{V}^2
    +\frac{L-C}2 \|w\|_{V'}^2
    +\frac12 \calI_1\\
   &\le \frac14 \|w\|_V^2 + C_L\|u_t\|_V^2 \|w\|_V^2.
   \label{Zelik24}
 \end{align}
 Finally, choosing $L$ so large that
 \beeq{Lsolarge}
     L\ge 2C \quand
   \frac{L}2\|w\|_{V'}^2 +\calI_1
    +\frac14 \|w\|_{V}^2 \ge 0,
 \end{equation}
 rewriting \eqref{Zelik24} (with obvious notation) as
 \beeq{newm}
   \ddt \calY + \kappa \calY
    \le m\calY,
 \end{equation}
 where
 \beeq{newm2}
   m:= C_L\|u_t\|_V^2 \in L^1(0,\infty)
 \end{equation}
 thanks to \eqref{relZelik1},
 the comparison principle for ODEs readily gives
 \eqref{relZelik2}.
\end{proof}
Note now that, comparing \eqref{dissin} and \eqref{relZelik2},
there follows in particular
\beeq{dissv}
  \int_t^\infty \|v_t(s)\|_{V'}^2\,\dis
   \le Q(\|U_0\|_2)e^{-\kappa t}
     +Q(\|g\|_{V'}).
\end{equation}
Thus, we can apply to \eqref{CHcomp} the procedure
used in Theorem~\ref{teosemi} to prove $\calV_2$-dissipativity.
Of course, the ``source'' term $LA^{-1}u_t$ in the \rhs\
of the differentiated equation
is easily controlled thanks to \eqref{dissin}.
Using also Lemma~\ref{lemmaR2} we then get
the estimate
\beeq{relZelik2b}
  \|V(t)\|_2^2+\|v_{tt}(t)\|^2_{V'}
   \le Q(\|U_0\|_2)e^{-t} + Q(\|g\|_{V'}).
\end{equation}
%
%
As a next step, we prove that the component $V$
of the solution is compact in $\calV_2$
and, more precisely, bounded in $\calV_3$. From this
point on, the further regularity \eqref{regogplus}
is needed.
\bele\label{lemmaZelik3}
 Let the assumptions of\/
 {\rm Theorem~\ref{teoattrastrong}} hold.
 Then we have
 \beeq{relZelik3}
   \|V(t)\|_3
   \le Q(\|U_0\|_2)e^{-\kappa t}
       +Q(\|g\|).
 \end{equation}
 %
\enle
\begin{proof}
 We differentiate \eqref{CHcomp} in time and test
 the result by $A(v_{tt}+\delta v_t)$ for small
 $\delta>0$. We do not give all the details, but
 just see how the nonlinear terms
 are controlled. Actually, performing some
 calculation and using
 \eqref{relZelik3} and interpolation,
 we get
 \bealo
   \big(f'(v)v_t,Av_{tt}\big)
    & =\frac12\ddt\io f'(v)|\nabla v_t|^2
     -\frac12\io f''(v)v_t|\nabla v_t|^2
     -\big( \dive (f''(v) v_t \nabla v), v_{tt}\big)\\
  \label{Zelik31}
    & \ge \frac12\ddt\io f'(v)|\nabla v_t|^2
    -\frac{\delta}4\|v_{t}\|_{D(A)}^2
    -\frac14\|v_{tt}\|^2
  -C_\delta\|v_t\|_V^2,
 \end{align}
 and, analogously,
 \beeq{Zelik31bis}
   \delta \big(f'(v)v_t,Av_{t}\big)
     \ge \delta \io f'(v)|\nabla v_t|^2
       -C\|v_t\|_V^2,
 \end{equation}
 where in both formulas $C$ (or $C_\delta$)
 is a monotone function of $\|V(t)\|_2$
 (and, more precisely, it depends at most
 polynomially on it).
 Thus, noting that the \rhs\ term
 $L(u_t,v_{tt}+\delta v_t)$ can be estimated
 in a standard way, one arrives
 at an expression of the form
 \beeq{Zelik31ter}
   \ddt \calY_3 + \kappa \calY_3
    \le Q(\|V(t)\|_2)
      + Q(\|U(t)\|_2)
    \le Q(\|U_0\|_2) e^{-\kappa t}
     + Q (\|g\|),
 \end{equation}
 where \eqref{dissifo} and \eqref{relZelik2b}
 have been used in deducing the latter inequality,
 and the functional $\calY_3$ (upon possibly
 taking a larger $L$) satisfies
 \beeq{Zelik31quater}
   c_L \|V_t\|_1^2 \le
     \calY_3 \le C_L \|V_t\|_1^2,
 \end{equation}
 where only $C_L$ depends on the radius of the
 absorbing set.
 Noting now that, by standard elliptic regularity
 results applied to \eqref{CHcomp}, we have
 \beeq{Zelik31quinquies}
   \|V(t)\|_3^2 \le C \big( \|V_t(t)\|_1^2
    + \|V(t)\|_2^2 \big),
 \end{equation}
 relation \eqref{relZelik3} comes then as an easy
 consequence of \eqref{Zelik31ter}. The lemma is proved.
\end{proof}
Finally, we show that $W$ is exponentially decaying
in $\calV_2$. Of course, this fact, together
with \eqref{relZelik3}, will give the desired property
of the decomposition \eqref{CHcomp}--\eqref{CHdecay}
and conclude the proof of Theorem~\ref{teoattrastrong}.
\bele\label{lemmaZelik4}
 Let the assumptions of\/
 {\rm Theorem~\ref{teoattrastrong}} hold.
 Then we have
 \beeq{relZelik4}
   \|W(t)\|_2
    \le Q(\|U_0\|_2) e^{-\kappa t}.
 \end{equation}
\enle
\begin{proof}
 We differentiate \eqref{CHdecay} in time and test
 the result by $w_{tt}+\delta w_t$ for small
 $\delta>0$. Still, the procedure is standard,
 but for the estimation of the nonlinear terms depending
 on $f$. Namely, we obtain on the \lhs
 \beeq{Zelik41}
   \big((f(u)-f(u-w))_t,w_{tt}+\delta w_t\big).
 \end{equation}
 Thus, defining
 \beeq{Zelik42}
   l=l(u,w):=\int_0^1 f'(su+(1-s)(u-w))\,\dis\ge -\lambda
 \end{equation}
 so that $f(u)-f(u-w)=lw$, we clearly have
 \beeq{Zelik43}
   \big((f(u)-f(u-w))_t,w_{tt}+\delta w_t\big)
    =\big(l_t w+l w_t,w_{tt}+\delta w_t\big)
 \end{equation}
 and
 \beeq{Zelik44}
  \big|\big(l_t w,w_{tt}+\delta w_t\big)\big|
   \le \|l_tw\|_V \|w_{tt}+\delta w_t\|_{V'}
   \le c\big(\|l_t\|_V\|w\|_{D(A)}\big)
   \|w_{tt}+\delta w_t\|_{V'}.
 \end{equation}
 Now, let us notice that,
 by \eqref{relZelik2b} and the analogue for $U$ coming
 from Theorem~\ref{teosemi} and Lemma~\ref{lemmaR2},
 \beeq{relZelik2c}
   \|W(t)\|_2^2 \le Q(\|U_0\|_2)e^{-t} +
    Q(\|g\|_{V'}).
 \end{equation}
 In particular, $\|l_t\|_V\le C$, with $C$ possibly
 depending on $U_0$, but independent of time.
 More precisely, using \eqref{relZelik2} and interpolation,
 we get, for all $\nu>0$,
 \beeq{relZelik2d}
    \|W(t)\|_{2-\nu}^2 \le Q(\|U_0\|_2)e^{-\kappa t},
 \end{equation}
 $\kappa$ depending here on $\nu$. Consequently
 (take $\nu=1$), we can control
 the \rhs\ of \eqref{Zelik44} so that
 \beeq{Zelik44.2}
  \big|\big(l_t w,w_{tt}+\delta w_t\big)\big|
    \le Q(\|U_0\|_2)e^{-\kappa t}
   +\frac14 \|w_{tt}+\delta w_t\|^2_{V'}
 \end{equation}
 and the latter term can be moved to the \lhs\ and
 estimated directly. Finally, coming back to the
 remaining term in \eqref{Zelik43}, we get
 \beeq{Zelik51}
   \big(l w_t,w_{tt}+\delta w_t\big)
    =\frac12\ddt \big(l,w_t^2\big)
      +\delta \big(l,w_t^2\big)
      -\frac12 \big(l_t,w_t^2\big),
 \end{equation}
 and the first two summands on the \rhs\ are controlled once
 more thanks to \eqref{f2}, while the third is estimated
 for small $\nu>0$ by
 \beeq{Zelik52}
   -\frac12 \big(l_t,w_t^2\big)
    \le c\|l_t\|_V \|w_t\|_{D\big(A^{\frac{1-\nu}2}\big)}^2
    \le C \|w_t\|_{D\big(A^{\frac{1-\nu}2}\big)}^2
    \le Q(\|U_0\|_2)e^{-\kappa t},
 \end{equation}
 thanks to \eqref{relZelik2d}.
 Thus, all the nonlinear terms
 are either (essentially) positive, or exponentially
 decaying. Then, \eqref{relZelik4} is proved,
 which concludes the proof of Lemma~\ref{lemmaZelik4}
 and of Theorem~\ref{teoattrastrong}.
\end{proof}

\section{Exponential attractors}
\label{secexpo}

This section is devoted to the proof of
existence of an exponential attractor for the
semiflow $\calS$ consisting of
the $\calV_2$-solutions to Problem~(P).
More precisely, we will prove
%
%
\bete\label{espattr}
 Assume\/ \eqref{f1}-\eqref{f3} and \eqref{regogplus}.
 Then, the semiflow $\calS$
 admits an exponential attractor $\calM_2$.
 Namely, $\calM_2$ is a positively invariant, compact subset
 of $\calV_2$ with finite fractal dimension
 with respect to the $\calV_2$-metric and bounded in $\calV_3$,
 such that, for any bounded $B\subset \calV_2$,
 there exist $C_B>0$ and $\kappa_B>0$ such that
 \beeq{expattr}
   \dist_2(S(t)B,\calM_2)\le C_B e^{-\kappa_B t},
 \end{equation}
 where $\dist_2$ denotes the\/ {\rm Hausdorff semidistance}
 of sets with respect to the $\calV_2$-metric.
\ente
Before proving the theorem, we need a couple
of preparatory lemmas.
\bele\label{prepZel1}
 Under the assumptions of\/ {\rm Theorem~\ref{espattr}},
 there exists a set $\calC_3$ bounded in $\calV_3$ which
 exponentially attracts any bounded set of $\calV_2$
 with respect to the $\calV_2$-metric.
\enle
\begin{proof}
 It is a simple consequence of the decomposition
 made in Section~\ref{secV2b}. More in detail,
 it follows from relations \eqref{relZelik3} and
 \eqref{relZelik4}.
\end{proof}
\bele\label{preplem2}
 Under the assumptions of\/ {\rm Theorem~\ref{espattr}},
 there exists a set $\calB_3$, bounded in $\calV_3$
 and positively invariant, which absorbs $\calC_3$
 and, consequently, exponentially attracts
 any bounded set of $\calV_2$ with respect to the
 distance of $\calV_2$.
\enle
\begin{proof}
 To prove the lemma, we basically need a dissipative
 estimate in $\calV_3$. This can be obtained just by
 mimicking the proof of Lemma~\ref{lemmaZelik3}. Namely,
 one has to differentiate \eqref{CH} in time and
 test the result by $u_{tt}+\delta u_t$ for small $\delta>0$.
 This leads to an expression perfectly analogous to
 \eqref{Zelik31ter}, with $\calY_3$ still satisfying
 \eqref{Zelik31quater}, but with $V$ everywhere replaced by
 $U$. This entails existence of a positively invariant
 and $\calV_3$-bounded set $\calB_3$, which eventually absorbs
 any $\calV_3$-bounded set of data. Since this in
 particular happens for $\calC_3$, the lemma is proved.
\end{proof}
\noindent%
{\bf Proof of Theorem~\ref{espattr}.}~~%
Let us start by considering initial data lying
in the set $\calB_3$ constructed above. Notice also
that it is not restrictive to assume $\calB_3$ to be
weakly closed in $\calV_3$.
Let us then take a couple of
solutions $u_1,u_2$ to Problem (P)
whose initial data
$(u\zzu,\uuu),(u\zzd,\uud)$ lie in $\calB_3$.
Since $\calB_3$ is positively invariant,
it is then clear that the functions
$t\mapsto (u_i(t),u_{i,t}(t))$, for $i=1,2$, take values in
$\calB_3$. By the $\calV_3$-analogue of \eqref{regouforte2},
which can be proved in a standard way,
we have that
\beeq{semprelimi}
  u_i\in C^2([0,\infty);H)
     \cap C^1([0,\infty);D(A))\cap C^0([0,\infty);D(A^2)),
\end{equation}
still for $i=1,2$.
Later on, the constants $c$ will be allowed to
depend on the choice of the initial
datum in $\calB_3$.

Let us now write equation \eqref{CH} for $u_1$ and $u_2$,
and then take the difference. This gives
\beeq{CHdiff}
  u\dtt+u_t+A(Au+f(u_1)-f(u_2))=0,
\end{equation}
where we have set $u:=u_1-u_2$. Let us then test
\eqref{CHdiff} by $A^{-1}(u_t+\delta u)$.
%
%
Setting
\beeq{defil}
  l=l(u_1,u_2)
   :=\int_0^1 f'(\tau u_1+(1-\tau) u_2)\,\ditau\ge -\lambda
\end{equation}
(cf.~\eqref{f2}) and writing $U:=(u,u_t)$,
standard manipulations lead us to the identity
\bealo
  & \ddt\Big[
   \frac12\|U\|^2_0
   +\frac\delta2\|u\|^2_{V'}
   +\frac12\io l(u_1,u_2)u^2
   +\delta \duav{u_t,A^{-1}u} \Big]\\
 \label{conto51new}
  & \mbox{}~~~~~
   +\delta\|u\|^2_V
   +(1-\delta)\|u_t\|^2_{V'}
   = -\delta\io l(u_1,u_2)u^2
   + \frac12 \io l_t(u_1,u_2)u^2.
\end{align}
Then, adding the inequality
\beeq{disuL}
  \ddt\big[ L\|u\|_{V'}^2 \big]
   \le L \big( \|u\|^2 + \|u_t\|_{D(A^{-1})}^2 \big),
\end{equation}
for $L$ large enough, using \eqref{f2}
(this also permits
to control the first term on the \rhs\
of \eqref{conto51new}),
and noting that, by standard use of
embeddings and interpolation,
\beeq{stimal}
  \frac12 \io l_t(u_1,u_2)u^2
   \le c \io \big (1+|u_{1,t}|+|u_{2,t}|\big)u^2
   \le \frac\delta2 \|u\|_V^2 + c_\delta \|u\|^2,
\end{equation}
relation \eqref{conto51new} takes the form
\beeq{conto53}
  \ddt\calY + \calZ
   \le c_{\delta,L}\big(\|u\|^2+\|u_t\|_{D(A^{-1})}^2\big)
    =c_{\delta,L}\big\|U\|_{-1}^2,
\end{equation}
with obvious
meaning of $\calY$ and $\calZ$.

Moreover, taking $\delta$ small enough and $L$ large
enough (the latter depending in particular
on~$\lambda$ in \eqref{f2}), it is clear
that, for some $c_3$, $c_4$, $\kappa$
also depending on the $\calV_3$-radius of $\calB_3$,
\beeq{Ysugiunew}
  c_3\|U\|_0^2
   \le \calY
   \le c_4\|U\|_0^2,
  \qquad \calZ \ge \kappa \calY,
\end{equation}
whence, taking $\ell>0$ and integrating
\eqref{conto53} from $\tau\in[0,\ell]$
to $2\ell$, we get
\beeq{conto54}
  \calY(2\ell)
   +\kappa \int_\tau^{2\ell}\calY(s)\,\dis
   \le \calY(\tau)
   +c_5\int_\tau^{2\ell}\|U(s)\|_{-1}^2\,\dis.
\end{equation}
(the $c_{\delta,L}$ in \eqref{conto53} has been
noted as $c_5$ for later convenience).
\newcommand{\Yy}[1]{\|(u,u_t)(#1)\|_0^2}%
\newcommand{\calUl}{{\mathcal U}_\ell}%
\newcommand{\calVl}{{\mathcal Z}_\ell}%
Now, let us apply the following straightforward fact
(see, e.g., \cite[Lemma~3.2]{Pr1}):
\bele\label{daPrazak}
 Let $\calH$ be a Hilbert space and
 $\calW$ a Banach space such that
 $\calH$ is compactly embedded into
 $\calW$. Then, for any $\gamma>0$,
 there exist a finite-dimensional orthonormal
 projector $P:\calH\to \calH$ and a positive
 constant $K$, both depending on $\gamma$ and
 such that, for all $w\in \calH$,
 \beeq{compastr}
   \|w\|^2_\calW \le \gamma\|w\|_\calH^2
    +K\|Pw\|_\calH^2.~~\dimbox
 \end{equation}
\enle
We apply Lemma~\ref{daPrazak} with
$\calH=\calV_0$ and $\calW=\calV_{-1}$.
Thus, with the notation above, we have
in particular
\beeq{comp}
  \|U\|_{-1}^2
   \le \gamma\|U\|_0^2
   +K\|PU\|_0^2.
\end{equation}
Let us now introduce the set of $\ell$-trajectories
associated with $\calV_0$-solutions of Problem~(P) as
\beeq{deficalul}
  \calUl:=\big\{
   U=(u,u_t)\in C^0([0,\ell],\calV_0):~U
   \text{ solves (P) on }\,[0,\ell]\big\}.
\end{equation}
The set $\calUl$ is endowed with the metric
of $L^2(0,\ell;\calV_0)$. For brevity,
we will write
\beeq{norcalul}
  \|U\|_\ell \text{~~in place of~~}
   \|U\|_{L^2(0,\ell;\calV_0)}.
\end{equation}
Note that, in general, $\calUl$ is not complete
with respect to~the chosen metric. However,
thanks to \eqref{conticompawweak}
of Theorem~\ref{teosemidebo}, it
is not difficult to see that, if
$\{U_n\}\subset \calUl$ is a sequence
such that $\{U_n(0)\}$ is {\sl bounded}\/
in $\calV_0$, and $U_n$ tends in
$L^2(0,\ell;\calV_0)$ to some function
$U$, then still
it is $U\in \calUl$. This is in fact our case
since we can restrict ourselves to
the subset $\calVl$ of the
elements of $\calUl$ whose initial
values lie in 
$\calB_3$. Actually, being $\calB_3$ {\sl weakly}\/
closed in $\calV_3$, it is easy to prove
that $\calVl$ is a complete
metric space with respect
to the $L^2(0,\ell;\calV_0)$-metric.
We can then define the shift operator
\beeq{defiL}
  \calL=\calL_\ell:\calVl\to\calVl, \qquad
   \calL(U)(\cdot):=U(\ell+\cdot).
\end{equation}
Integrating now \eqref{conto54} with respect to~$\tau$
from $0$ to $\ell$ and using \eqref{Ysugiunew}, \eqref{comp},
we infer
\bealo
  & c_3\ell\|U(2\ell)\|_0^2
   +c_3\kappa\ell\int_\ell^{2\ell}
    \|U(s)\|_0^2\,\dis\\
 \label{conto55}
  & \mbox{}~~~~~
   \le c_4\int_0^\ell\|U(\tau)\|_0^2\,\ditau
   +c_5\gamma \ell\int_0^{2\ell}\|U(s)\|_0^2\,\dis
   +c_5K\ell\int_0^{2\ell} \|PU(s)\|_0^2\,\dis.
\end{align}
Then, dividing \eqref{conto55} by $\ell$
and using the notation \eqref{defiL} we obtain
\beeq{conto56}
  (c_3\kappa-c_5\gamma) \|\calL U\|_\ell^2
   \le \Big(\frac{c_4}\ell+c_5\gamma\Big)\|U\|_\ell^2
   +c_5 K\big(\|PU\|_\ell^2+\|P\calL U\|_\ell^2\big).
\end{equation}
Now, let us choose in turn $\gamma$ and $\ell$
such that
\beeq{gammaell}
  c_5\gamma \le \frac{c_3\kappa}{17}, \qquad
   \frac{c_4}{\ell}\le \frac{c_3\kappa}{17}.
\end{equation}
Thus, for some $K'$ depending on all other
constants, \eqref{conto56} gives,
for all $u\in\calVl$,
\beeq{conto57}
  \|\calL U\|_\ell^2
   \le \frac18\|U\|_\ell^2
   +K'\big(\|PU\|_\ell^2+\|P\calL U\|_\ell^2\big).
\end{equation}
Consequently, the semiflow $\calS$ associated with
Problem (P) enjoys the {\sl generalized squeezing
property}\/ introduced in \cite[Def.~3.1]{Pr1}
on the set $\calB_3$. Recalling
\cite[Lemma~2.2]{Pr2}, we infer that the
{\sl discrete}\/ dynamical system on $\calVl$
generated by $\calL$ admits an exponential attractor
$\calM_{\discr}$, which is compact
and has finite fractal dimension in $\calV_0$
and exponentially attracts $\calB_3$ in the
$\calV_0$-metric.

To pass from $\calM_{\discr}$
to a {\sl regular}\/ exponential attractor
for the original semiflow, we proceed by noting
a number of facts:\\[2mm]
{\sl (a)~~}%
The evaluation map $e:\calVl\to \calV_0$ given by
$e:U\mapsto U(\ell)$ is Lipschitz continuous.
To prove this, one can, e.g., multiply \eqref{conto53}
by $t$ and integrate in time between
$0$ and $\ell$.\\[1mm]
{\sl (b)~~}%
The semigroup operator $S(t)$ is uniformly
Lipschitz continuous on $[0,\ell]$
with respect to the metric of $\calV_0$.
This is easily shown by integrating
once more \eqref{conto53} between
$0$ and an arbitrary $t\le \ell$ and using
Gronwall's lemma.\\[1mm]
{\sl (c)~~}%
For each solution $U\in\calVl$ and all
$0\le s\le t\le \ell$, by \eqref{regouforte}
there holds
\beeq{Holderintime}
  \|U(t)-U(s)\|_0^2
   \le \Big|\int_s^t \|U_t(\tau)\|_0\,\ditau\Big|^2
   \le c|t-s|^2.
\end{equation}
In other words, the $\ell$-trajectories lying
in $\calVl$ are uniformly Lipschitz continuous in time.\\[1mm]
{\sl (d)~~}%
Properties {\sl (a)}-{\sl (c)}\/ allow us to apply,
e.g., \cite[Thm.~2.6]{MP} to deduce
that there exists an invariant compact subset
$\calM_0\subset \calV_0$, of finite fractal dimension
with respect to the $\calV_0$-topology,
which exponentially attracts $\calB_3$
still with respect to the $\calV_0$-topology.
More precisely, since the elements of $\calVl$
take values in the $\calV_3$-bounded and
positively invariant set $\calB_3$,
setting $\calM_2:=\calM_0\cap \calB_3$,
we have that $\calM_2$ is bounded in $\calV_3$
and, by interpolation, it is compact
and has finite fractal dimension in
$\calV_2$ (in fact, in $\calV_s$ for any $s<3$).
Moreover, it exponentially attracts $\calB_3$
and, by interpolation, this happens even with respect to
the $\calV_2$-metric;\\[1mm]
{\sl (e)~~}%
Finally, we see that $\calM_2$ exponentially attracts
any set $B$ bounded in $\calV_2$. Actually, we
know from point~{\sl (d)} that
$\calM_2$ exponentially attracts $\calB_3$
and from Lemma~\ref{preplem2} that
$\calB_3$ exponentially attracts any such $B$.
Note that the exponential attraction
holds in both cases
with respect to the $\calV_2$-metric.
To conclude, we can thus apply the {\sl transitivity
property}\/ of exponential attraction introduced in
\cite[Thm.~5.1]{FGMZ}. To do this, we have to check
(cf.~\cite[(5.1)]{FGMZ}) that the semigroup operators $S(t)$
are uniformly Lipschitz continuous on bounded balls $B$ of
$\calV_2$, with the Lipschitz constant having the form
$c_6e^{c_7t}$, where $c_6$ and $c_7$ depend only on $B$. To
prove this fact, we can rewrite \eqref{conto31} for two
solutions $u^1$ and $u^2$ originating from initial data
$U^1_0=(u_0^1,u_1^1)$ and $U^2_0=(u_0^2,u_1^2)$, respectively.
We obtain (recall also \eqref{CH2})
$$
\|U(t)\|_2^2 +\itt\|\nabla U_t\|^2
\le\big\|U^1_0-U^2_0\big\|_2^2 +\itt\|\nabla (G^1-G^2)\|^2.
$$
Then, it is not difficult to recover the wanted estimate.
The proof is completed. \dimbox
%
%
%
%


      \section{Energy solutions}
       \label{secV0}


We finally consider the class of energy solutions. As we shall see, in this case
the dissipation integral \eqref{dissin} will not be used.

We start by establishing existence and uniqueness of solutions.
\bete\label{teoesidebo}
 Let us assume\/ \eqref{f1}-\eqref{f3}
 and \eqref{regog}, together with
 \beeq{regou0}
    (u_0,u_1)\in\calV_0.
 \end{equation}
 Then, there exists\/ {\rm one and only one} function
 \beeq{regou}
   u\in W^{1,\infty}(0,T;V')\cap L^\infty(0,T;V)
 \end{equation}
 which solves\/ {\rm Problem~(P)}.
\ente
\beos\label{remgrowth}
 It will be clear from the proof that the growth restriction
 \eqref{f3} is required only for uniqueness. Actually,
 existence in the class $\calV_0$ holds for
 any polynomial growth (cf., e.g., \cite[Thm.~1.1]{Ba2}).
\eddos
\noindent%
{\bf Proof of Theorem~\ref{teoesidebo}.}~~
Let us start with the proof of existence, which
follows closely \cite[Sec.~4]{S} and is reported
just for later convenience.
We let $\{\lambda_j\}\subset(0,\infty)$, $j\in\NN$, be the
sequence of eigenvalues of $A$, increasingly
ordered and with possible repetitions according to the
multiplicities. Correspondingly, we let $\{z_j\}$ be
a (complete) systems of eigenvectors, which is chosen
to be orthonormal in $H$ and orthogonal in $V$. We
set $Z_N:=\spa\{z_1,\dots,z_N\}$ and denote
by $P_N$ the orthogonal projector onto $Z_N$. Of course,
$P_N$ can be thought to act on any of the spaces
$D(A^s)$, $s\in\RR$. We then consider the Faedo-Galerkin
approximation of Problem~(P), i.e.,
\begin{align}\label{CHGal}
  & A^{-1}\big(v_{N,tt}+v_{N,t}\big)
    +A v_N+P_N f(v_N)=P_N A^{-1} g, \\
 \label{inizGal}
  & V_N|_{t=0}=V_{0,N}:=P_N(u_0,u_1),
\end{align}
where both relations are intended as
equalities in $Z_N$ and we have set,
for brevity, $V_N:=(v_N,v_{N,t})$.
It is easy to show that
Problem~\eqref{CHGal}-\eqref{inizGal}
admits one and only one solution, which satisfies
the energy estimate
(cf.~\eqref{conto11}) uniformly
with respect to~$N$. Then, standard compactness tools
and the growth restriction~\eqref{f3}
(note that at this level any polynomial growth
of $f$ would be admissible, cf.~Remark~\ref{remgrowth})
permit to take the limit of~\eqref{CHGal}-\eqref{inizGal}
as at least a subsequence of $N$ goes to $\infty$
and, as a consequence, get existence of one
$\calV_0$-solution $U=(u,u_t)$ to Problem~(P)
satisfying in particular~\eqref{regou}.

To get uniqueness, following the method developed in
\cite{Sed}, we will prove that,
as $U=(u,u_t)$ is {\sl any}\/ solution to~(P)
in the regularity class~\eqref{regou},
the {\sl whole sequence}\/
$V_N$ converges to $U$.
Of course, this entails uniqueness
of $U$. With this aim, we let
$U_N:=P_N U$ (i.e., $u_N:=P_N u$ and
$u_{N,t}:=P_N u_t$) and consider the projection
of equation \eqref{CH}.
%
%
%
Then, it is clear that the difference
$W_N=(w_N,w_{N,t}):=U_N-V_N$ satisfies
\begin{align}\label{CHGaldiff}
  & A^{-1}\big(w_{N,tt}+w_{N,t}\big)
    +A w_N+P_N\big(f(u_N)-f(v_N)\big)
    =P_N\big(f(u_N)-f(u)\big),\\[1mm]
 \label{inizGaldiff}
  & W_N|_{t=0}=0.
\end{align}
Testing \eqref{CHGaldiff} by $A^{-1}w_{N,t}$,
we readily get
\bealo
  & \frac12\ddt\|W_N\|_{-1}^2
    +\|w_{N,t}\|_{D(A^{-1})}^2\\
  \nonumber
   &= \duavg{P_N\big(f(v_N)-f(u_N)\big),A^{-1}w_{N,t}}
   +\duavg{P_N\big(f(u_N)-f(u)\big),A^{-1}w_{N,t}}\\
   &\le \|f(u_N)-f(v_N)\|\|w_{N,t}\|_{D(A^{-1})}
    +\|f(u_N)-f(u)\|_{V'}\|w_{N,t}\|_{V'}.
  \label{contoa1}
\end{align}
Let us then notice that,
by \eqref{f3} and for fixed but
arbitrary $\epsilon>0$, one has
\bealo
  \|f(u_N)-f(u)\|_{V'}
   & \le c\Big\| \int_0^1 f'(\tau u_N+(1-\tau)u)\,\ditau
        (u_N-u)\Big\|_{L^{1+\epsilon}(\Omega)}\\
 \label{sergey1}
   & \le C\|u_N-u\|
   \le C\lambda_N^{-1/2}\|u_N-u\|_V
   \le C \lambda_N^{-1/2}.
\end{align}
Here and below, the constant $C$ is allowed to
depend on the $L^\infty(0,T;\calV_0)$-norms of
$U$ and $V_N$ (of course, they are bounded independently
of $N$). Thus, using once more \eqref{f3} and
the Br\'ezis-Gallouet inequality (cf.~\eqref{interpo2}
or \eqref{interpo3}),
the remaining term in \eqref{contoa1}
can be controlled as
\bealo
  \|f(u_N)-f(v_N)\|
   & \le C\big(1+\|u_N\|^2_{L^\infty(\Omega)}
             +\|v_N\|^2_{L^\infty(\Omega)}\big)
     \|w_N\|\\
 \label{sergey2}
   & \le C\big(1+\log(1+\|u_N\|_{D(A)})
             + \log(1+\|v_N\|_{D(A)})\big) \|w_N\|
     \le C\log \lambda_N \|w_N\|.
\end{align}
Thus, using \eqref{sergey1} and \eqref{sergey2},
\eqref{contoa1} yields
\beeq{contoa2}
  \frac12\ddt\|W_N\|_{-1}^2
    +\|w_{N,t}\|_{D(A^{-1})}^2
   \le C\log \lambda_N \|w_N\|\|w_{N,t}\|_{D(A^{-1})}
   + C \lambda_N^{-1/2},
\end{equation}
whence in particular
\beeq{contoa3}
  \ddt\|W_N\|_{-1}^2
     \le C_1\log \lambda_N \|W_N\|_{-1}^2
   +C_2\lambda_N^{-1/2},
\end{equation}
and, by Gronwall's lemma,
\beeq{contoa4}
  \|W_N(t)\|_{-1}^2
   \le C_2 t \lambda_N^{C_1 t-\frac12},
\end{equation}
so that, taking, e.g., $t_*:= 1/4C_1$, we readily
obtain that, as (the whole sequence) $N\nearrow\infty$,
\beeq{contoa5}
  W_N \to 0 \quext{strongly in }\,
   L^\infty(0,t_*;\calV_{-1}),
\end{equation}
whence, since we already know that $U_N\to U$,
we obtain $V_N\to U$ by comparison.
Finally, restarting the procedure from
the time $t_*$ (and noting that the value
of the ``new'' $t_*$ does not change since
the functions $U_N,V_N$ stay bounded in $\calV_0$
uniformly in time), we deduce uniqueness
on the whole of $(0,T)$,
which concludes the proof. \dimbox

\vspace{2mm}
\noindent%
Therefore, the energy solutions constitute a new
semiflow~$\calS_0$. The following analogue
of Theorem~\ref{teosemi} establishes some
properties of $\calS_0$ and of the associated
semigroup operator $S_0$.
\bete\label{teosemidebo}
 Let the assumptions of\/ {\rm Theorem~\ref{teoesidebo}}
 hold. Then, the semiflow~$\calS_0$ is\/ {\rm uniformly
 dissipative}. Namely, there exists a constant $R_0$
 independent of the initial data such that, for all
 bounded $B\subset \calV_0$, there exists $T_B\ge 0$
 such that $\|S_0(t)b\|_0\le R_0$, for all
 $b\in B$ and $t\ge T_B$.
 Moreover, any $u\in \calS_0$ satisfies the additional
 time continuity property
 \beeq{regouforte2weak}
   u\in C^2([0,T];D(A^{-3/2}))\cap C^1([0,T];V')\cap C^0([0,T];V)
 \end{equation}
 as well as the\/ {\rm energy equality}
 \beeq{energyeq}
   \calE(u,u_t)(t)-\calE(u,u_t)(s)
    =-\int_s^t\|u_t(r)\|_{V'}^2\,\dir
    \qquad\perogni s,t\in [0,T].
 \end{equation}
 Finally, given a sequence of initial data
 $\{(u\znn,u\unn)\}\subset \calV_0$ suitably tending to
 some $(u_0,u_1)\in \calV_0$,
 and denoting by $u_n,u$ the solutions
 emanating from $(u\znn,u\unn),(u_0,u_1)$, respectively,
 we have that
 \beal{conticompawweak}
   (u\znn,u\unn)\to (u_0,u_1) \text{~~weakly in }\/\calV_0
    & ~\Rightarrow~ (u_n,u_{n,t}) \to (u,u_{t})
    \text{~~weakly star in }\,L^\infty(0,T;\calV_0),\\
 \label{conticompaweak}
   (u\znn,u\unn)\to (u_0,u_1) \text{~~strongly in }\/\calV_0
    & ~\Rightarrow~ (u_n,u_{n,t}) \to (u,u_{t})
    \text{~~strongly in }\,C^0([0,T];\calV_0),
 \end{align}
 for any fixed $T\ge 0$.
\ente
\noindent%
\begin{proof}
We start by showing \eqref{regouforte2weak}.
Let $u$ be the $\calV_0$-solution to (P) and
set $v:=e^{t/2}u$ so that
\beeq{derivv}
  v_t=e^{t/2}u_t+\frac12e^{t/2}u
     =e^{t/2}u_t+\frac v2, \qquad
  v_{tt}=e^{t/2}u_{tt}+\frac12e^{t/2}u_t+\frac{v_t}2
        =e^{t/2}u_{tt}+v_t-\frac v4.
\end{equation}
Let us now (formally)
multiply \eqref{CH} by $e^{t/2}A^{-1}v_t$.
After some calculations we obtain
\beeq{contob1}
  \ddt\calY(t)=\Phi(t):= e^t\io\big(2F(u)-f(u)u)
   -e^t\duav{g,A^{-1}u},
\end{equation}
where we have set
\beeq{contob2}
  \calY(t):= \|v_t\|_{V'}^2
   +\|v\|_{V}^2
   -\frac14\|v\|_{V'}^2
   +2e^t\io F(u)
   -2e^t\duav{g,A^{-1}u}.
\end{equation}
Of course, \eqref{contob1} could make no sense, because
$e^{t/2}A^{-1}v_t$ is not smooth enough to be
used as a test function. However,
if we let $u_n$ be a class of $\calV_2$-solutions
suitably approximating $u$ and define $v_n$ accordingly,
it is clear that then \eqref{contob1} holds
at least for $v_n$.

More precisely, we can suppose that, given
$s,t\in[0,T]$, with $s<t$, there holds
at least
\beeq{unudebo}
  u_n\to u \quext{weakly star in }\,
   W^{1,\infty}(s,t;V') \cap L^\infty(s,t;V).
\end{equation}
Thus, by \eqref{f3} and
Lebesgue's dominated convergence theorem,
it readily follows that
\beeq{phiconve}
  \Phi_n\to \Phi \quext{strongly in }\,
    L^1(s,t)
\end{equation}
(with obvious meaning of $\Phi_n$,
cf.~\eqref{contob1}).
To proceed, we additionally assume that
\beeq{unuS}
  (u_n,u_{n,t})(s)\to (u,u_t)(s) \quext{strongly in }\,
   V\times V'.
\end{equation}
Clearly, this can be done as
we consider $s$ as the initial time
and choose the approximation
$u_n$ accordingly. Then, integrating
\eqref{contob1} (written for $u_n$) over $(s,t)$,
taking the supremum limit
as $n\nearrow\infty$, and using
\eqref{unudebo} and the trivial fact that
$\calY$ (seen as a functional of the couple
$(v,v_t)$) is sequentially weakly lower semicontinuous
in $\calV_0$, one gets that for the limit solution
$u$ there holds
\beeq{disuene}
  \calY(t)\le \calY(s)
   +\int_s^t\Phi(r)\,\dir.
\end{equation}
To prove the converse inequality, one
simply repeats the procedure by
considering $t$ as the initial time and noting that,
due its hyperbolic nature, (P) is solvable backward
in time (of course, this prevents dissipation,
but still there is global boundedness
in the energy norm).
In particular, the approximation $u_n$ can be still
chosen to fulfill \eqref{unudebo}, while
in place of \eqref{unuS} we can ask that
\beeq{unuT}
  (u_n,u_{n,t})(t)\to (u,u_t)(t) \quext{strongly in }\,
   V\times V'.
\end{equation}
Thus, we finally get the {\sl equal sign}\/ in \eqref{disuene},
which, due to arbitrariness of $s,t$, readily implies
that $\calY$ (written for the limit solution $u$ and regarded
as a function of time) is absolutely continuous over
$[0,T]$.
To conclude, we observe that from \eqref{regou} we know
that
\beeq{uCw}
  v\in C^1_w([0,T];V')\cap C^0_w([0,T];V),
\end{equation}
so that the latter three summands
in $\calY$ (cf.~\eqref{contob2})
are strongly continuous in time. Then, by comparison,
the function $t\mapsto \|v_t(t)\|_{V'}^2+\|v(t)\|_{V}^2$
is also continuous in time. This fact, joint with
\eqref{uCw}, immediately gives \eqref{regouforte2weak}
(as before, the continuity of $u_{tt}$ can
be shown by a further comparison of terms).

We now show that $u$ satisfies \eqref{energyeq},
i.e., the energy equality for the original
energy~$\calE$. We give just the highlights of the
argument and leave the details to the reader.
Setting $\calY_0:=e^{-t}\calY$ and $\Phi_0:=e^{-t}\Phi$,
and noting that $\calY$, and hence $\calY_0$,
is absolutely continuous in time, we infer
from \eqref{contob1}
\beeq{intermeY}
  \calY_0'+\calY_0=\Phi_0 \quext{a.e.~in }\,(0,T).
\end{equation}
Next, let us integrate \eqref{intermeY} between $s$ and $t$,
and use \eqref{derivv}-\eqref{contob2} (in particular,
$\calY_0$ and $\Phi_0$ have to be rewritten in terms of
$u,u_t$ rather than $v,v_t$). Performing standard
manipulations and subtracting from the resulting formula
the outcome of \eqref{CH} tested by $A^{-1}u$, one
then gets \eqref{energyeq} by
simple computations. It is maybe worth
pointing out that \eqref{regouforte2weak},
or even \eqref{regou}, is sufficient to test
\eqref{CH} by $A^{-1}u$ (while it does {\sl not}\/
permit to test it by $A^{-1}u_t$), so the latter
argument is rigorous.

As we know the energy equality to hold,
the dissipativity of $S_0$ can be standardly
obtained by writing \eqref{energyeq} in the differential
form (which is possible almost everywhere by absolute continuity
of $\calE$) and adding the result of
\eqref{CH} tested by $\delta A^{-1} u$ for small
$\delta>0$. We leave once more the details to
the reader.

Finally, we have to prove \eqref{conticompawweak} and
\eqref{conticompaweak}. The first is standard. To show
the latter, we proceed along the lines of
\cite[Proof of Thm.~3.6]{Ba2}. Namely,
setting as usual $U_n:=(u_n,u_{n,t})$, we
let by contradiction $\epsilon>0$ and
$\{t_n\}\subset [0,T]$ such that
$\|U_n(t_n)-U(t_n)\|_0\ge \epsilon$
at least for a subsequence (here we shall
not relabel subsequences).
We can also assume $t_n\to t$ for some
$t\in[0,T]$. Note that then, by weak
convergence, $U_n(t_n)\to U(t)$ weakly in $\calV_0$.

Now, we claim it suffices to prove that
$|\calE(U_n(t_n))-\calE(U(t))|\to 0$.
Indeed, since the
other terms in $\calE$ have a lower order,
this would entail
\beeq{Ball0}
  \big| \|U_n(t_n)\|_0^2 - \|U(t)\|_0^2 \big| \to 0
\end{equation}
and consequently that $U_n(t_n)$ goes to $U(t)$
{\sl strongly}\/ in $\calV_0$. Now, by the
energy equality it is
\beeq{Ball0.5}
  \big|\calE(U(t_n))-\calE(U(t))\big|
   \le \Big| \int_{t_n}^t \|u_t\|_{V'}^2 \Big|\to 0,
\end{equation}
so that $U$ is strongly continuous in time and
in particular $U(t_n)\to U(t)$ strongly in $\calV_0$.
Thus, we would end up with
\beeq{Ball0.6}
  \big\|U_n(t_n) - U(t_n)\big\|_0 \to 0,
\end{equation}
a contradiction.

Let us now prove the claim. First,
by sequential weak lower semicontinuity
of $\calE$, we have
\beeq{Ball1}
  \calE(U(t)) \le \liminf_{n\nearrow\infty}
   \calE(U_n(t_n)).
\end{equation}
Conversely, by the energy equality, we obtain
\bealo
  \limsup_{n\nearrow\infty} \calE(U_n(t_n))
    & = \lim_{n\nearrow\infty} \calE(U_{0,n})
     - \liminf_{n\nearrow\infty} \int_0^{t_n} \|u_{n,t}\|_{V'}^2\\
 \no
   & = \calE(U_0)
     - \liminf_{n\nearrow\infty} \int_0^t \|u_{n,t}\|_{V'}^2
      - \lim_{n\nearrow\infty} \int_t^{t_n}\|u_{n,t}\|_{V'}^2\\
 \label{Ball2}
   & \le \calE(U_0)
    - \int_0^t \|u_t\|_{V'}^2
     - \lim_{n\nearrow\infty} \int_t^{t_n}\|u_{n,t}\|_{V'}^2
     = \calE(U(t)),
\end{align}
provided that the last integral on the second and
third row does go to $0$, which is true for
$t_n\to t$ and
\beeq{Ball3}
   \limsup_{n\nearrow\infty}\|u_{n,t}\|^2_{L^2(0,T;V')}
    \le \calE(U_0)-\liminf_{n\nearrow\infty}\calE(U_n(T))
    \le \calE(U_0)-\calE(U(T))
    =\|u_t\|^2_{L^2(0,T;V')},
\end{equation}
i.e., $u_{n,t}$ goes to $u_t$ {\sl strongly}\/
in $L^2(0,T;V')$.
The proof is thus complete.
\end{proof}
\vspace{2mm}
\noindent%
Still paralleling the $\calV_2$-case, we finally have the
\bete\label{teoasydebo}
 Let the assumptions of\/ {\rm Theorem~\ref{teoesidebo}}
 hold. Then, the semiflow $\calS_0$ associated
 with~{\rm (P)} is\/ {\rm asymptotically compact}.
 Thus $\calS_0$ possesses the global attractor $\calA_0$.
\ente
\begin{proof}
The proof is analogous to that of
Theorem~\ref{teoasy} and, in fact,
even technically simpler. Indeed, having the
$\calV_0$-energy equality \eqref{energyeq}
at our disposal, we can still implement
Ball's energy method \cite[Sec.~4]{Ba2},
which leads us to show the asymptotic compactness of $S_0$
in the phase space $\calV_0$. Then we can conclude as in Section
\ref{secV2}.
\end{proof}
\noindent \beos\label{weakf} It is not difficult to check that
the results of this Section still hold when $f\in C^1({\mathbb
R};{\mathbb R})$ only (compare with \eqref{f1}). \eddos
\noindent \beos\label{attr0-stru} The semiflow $\calS_0$ is
generated by a gradient system (see \eqref{defiE}). Hence
$\calA_0$ coincides with the unstable manifold of the set of
equilibria. We recall that this set can have a very complicated
structure (see, e.g., \cite{BF,MMW,WW1,WW2,WW3}). Moreover, it
is clear that $\calA_2 \subset \calA_0$. However, the converse
is far less trivial (see \cite{GSZ2}). \eddos
\noindent
\beos\label{LS}
 Any energy solution given by Theorem~\ref{teoesidebo} converges to a unique
 equilibrium, provided that $f$ is real analytic.
 This fact can be proven, using Theorem~\ref{teoasydebo} and
 the {\L}ojasiewicz-Simon inequality,
 arguing as in \cite{GPS}. A convergence rate
 estimate can also be obtained.
\eddos
\beos\label{ancheV1}
 Well-posedness of \eqref{CH} can be also shown in
 the space $\calV_1$, i.e., for weak solutions. Actually, the uniqueness
 part of Theorem~\ref{teoesiforte} still holds
 with no change in the proof. Concerning existence,
 the main bound should be obtained
 by testing \eqref{CH} by $u_t+\beta u$,
 for small $\beta>0$ (compare this
 with \eqref{conto21} below). However, this estimate
 does not seem to have a dissipative character
 since one apparently cannot take advantage of
 the dissipation integral \eqref{dissin}.
 In this sense, the $\calV_1$-theory for
 Problem~(P) seems less complete than the
 $\calV_2$ and $\calV_0$ theories discussed above.
 We also notice that, both in the $\calV_2$
 and in the $\calV_1$ setting, it seems possible
 to obtain global well posedness for nonlinearities
 $f$ having up to a fourth order growth
 (rather than a cubic growth as stated
 in~\eqref{f3}). However, in this case we
 can no longer prove dissipativity,
 even for quasi-strong solutions.
\eddos


\bigskip
\noindent
\textbf{Acknowledgments.} The first two authors have been partially supported by
the Italian PRIN Research Project 2006 {\sl ``Problemi a frontiera libera,
transizioni di fase e modelli di isteresi''}.



\end{document}